\documentclass[a4paper,12pt]{amsart}

\usepackage{amssymb}
\usepackage{latexsym}
\usepackage{amsmath}
\usepackage{amscd}
\usepackage{times}

\title[The second rational cohomology 
of the Johnson kernel]{The second Johnson 
homomorphism and the second 
rational cohomology of the Johnson kernel}

\author{Takuya Sakasai}
\address{Graduate School of Mathematical Sciences, 
the University of Tokyo, 3-8-1 Komaba, Meguro-ku, 
Tokyo 153-8914, Japan}
\email{sakasai@ms.u-tokyo.ac.jp}


\newtheorem{thm}{Theorem}[section]
\newtheorem{prop}[thm]{Proposition}
\newtheorem{lem}[thm]{Lemma}

\newtheorem{tthm}{Theorem \ref{mainSg}\hspace{-4pt}}

\newtheorem{llem}{Lemma \ref{irrSg1}\hspace{-4pt}}

\theoremstyle{definition}

\newtheorem{remark}[thm]{Remark}

\begin{document}

\newcommand{\Mg}{\mathcal{M}_g}
\newcommand{\Mgp}{\mathcal{M}_{g,\ast}}
\newcommand{\Mgb}{\mathcal{M}_{g,1}}
\newcommand{\Ig}{\mathcal{I}_g}
\newcommand{\Igp}{\mathcal{I}_{g,\ast}}
\newcommand{\Igb}{\mathcal{I}_{g,1}}
\newcommand{\Kg}{\mathcal{K}_g}
\newcommand{\Kgp}{\mathcal{K}_{g,\ast}}
\newcommand{\Kgb}{\mathcal{K}_{g,1}}
\newcommand{\Sg}{\Sigma_g}
\newcommand{\Sgb}{\Sigma_{g,1}}
\newcommand{\Symp}[1]{Sp(2g,\mathbb{#1})}
\newcommand{\symp}[1]{\mathfrak{sp}(2g,\mathbb{#1})}
\newcommand{\DHg}{\mathfrak{h}_g}
\newcommand{\DHgb}{\mathfrak{h}_{g,1}}
\newcommand{\DHgp}{\mathfrak{h}_{g,\ast}}
\newcommand{\Lg}{\mathcal{L}_g}
\newcommand{\Lgb}{\mathcal{L}_{g,1}}
\newcommand{\At}[1]{\mathcal{A}_{#1}^t (H)}
\newcommand{\Hq}{H_{\mathbb{Q}}}
\newcommand{\Ker}{\mathop{\mathrm{Ker}}\nolimits}
\newcommand{\Hom}{\mathop{\mathrm{Hom}}\nolimits}

\newcommand{\Out}{\mathop{\mathrm{Out}}\nolimits}
\newcommand{\Aut}{\mathop{\mathrm{Aut}}\nolimits}
\renewcommand{\Im}{\mathop{\mathrm{Im}}\nolimits}
\newcommand{\Q}{\mathbb{Q}}
\newcommand{\Z}{\mathbb{Z}}

\begin{abstract}
The Johnson kernel is the subgroup of the mapping class 
group of a surface generated by Dehn twists along 
bounding simple closed curves, and has the second Johnson 
homomorphism as a free abelian quotient. 
In terms of the representation 
theory of the symplectic group, 
we give a complete description of cup products of 
two classes in the first rational cohomology of the Johnson kernel 
obtained by the rational dual of 
the second Johnson homomorphism. 
\end{abstract}

\renewcommand\baselinestretch{1.1}
\setlength{\baselineskip}{16pt}

\newcounter{fig}
\setcounter{fig}{0}

\maketitle

\section{Introduction}\label{sec:intro}
Let $\Sg$ be a closed oriented surface of genus $g \ge 2$ and 
let $\Mg$ be the mapping class group of $\Sg$, 
which is the group of isotopy classes of 
orientation preserving diffeomorphisms of $\Sg$. 
The Torelli group $\Ig$ is the subgroup of $\Mg$ consisting of 
all elements which act trivially on the first homology group 
$H:=H_1 (\Sg)$ of $\Sg$, and the group $\Kg$ is the subgroup of 
$\Ig$ generated 
by Dehn twists along bounding simple closed curves. 
Johnson \cite{jo1} showed that 
$\Kg$ coincides with the kernel of what is now called 
the first Johnson homomorphism 
$\tau_g (1)$ of $\Ig$ (see \cite{jo}), 
namely we have an exact sequence
\[\begin{CD}
1 @>>> \Kg @>>> \Ig @>{\tau_g (1)}>> \DHg (1)=(\wedge^3 H)/H @>>> 1, 
\end{CD}\]
and by this, $\Kg$ is called the {\it Johnson kernel} or 
the {\it Johnson subgroup}. 

The group $\Kg$ plays an important role in topology. For 
example, it has some relationships to 
the Casson invariant of homology 3-spheres and secondary characteristic 
classes of surface bundles as we see in Morita's papers \cite{mo9, mo6}. 
However, we still do not have enough information on $\Kg$. 
McCullough-Miller \cite{mm} showed that $\mathcal{K}_2=\mathcal{I}_2$ is 
not finitely generated, and Mess \cite{mess} showed that 
it is a free group of infinite rank. 
Recently, Biss-Farb \cite{bif} showed that $\Kg$ is not 
finitely generated for all $g \ge 2$. Note that the determination of 
the abelianization of $\Kg$ is still open. 

Both of the roles of $\Kg$ mentioned above can be interpreted as properties 
of some elements in the rational cohomology $H^\ast (\Kg;\Q)$ of $\Kg$, to 
which we now pay our attention. 
By the fact that $\Kg$ is torsion-free and acts on 
the Teichm\"uller space 
properly discontinuously, 
we can see that $\Kg$ has finite cohomological dimension. On the 
other hand, Akita \cite{akita} showed that $H^\ast (\Kg;\Q)$ is an 
infinite-dimensional vector space. 

Our strategy to study $H^\ast (\Kg;\Q)$ is 
to use the second Johnson homomorphism 
\[\tau_g (2): \Kg \longrightarrow \DHg (2),\]
where $\DHg (2)$ is a certain free abelian group which can be 
described by using $H$ (see Section \ref{sec:preliminaries} 
for the notation used here). 
An important fact of $\tau_g (2)$ is that it is $\Mg$-equivariant, 
where $\Mg$ acts on $\Kg$ by conjugation and acts on $\DHg (2)$ 
through its action on $H$ known as the classical 
representation $\Mg \to \Symp{Z}$. Moreover, if we consider 
the induced map 
\[\tau_g (2)_{\ast}: H_1(\Kg;\Q) 
\longrightarrow H_1 (\DHg (2);\Q)=\DHg (2) \otimes \Q,\]
we can check that its image is preserved by the action of 
$\Symp{Q}$ extending that of $\Symp{Z}$ 
(see Asada-Nakamura cite[Lemma 2.2.8]{an}). 
In terms of the symplectic representation theory, 
$\DHg (2) \otimes \Q$ is the irreducible representation 
denoted by $[2^2]$ (see \cite{hain}, \cite{mo9}). 
By passing to 
the dual, we obtain an injection 
\begin{align*}
\tau_g (2)^\ast : [2^2] &\longrightarrow H^1(\Kg;\Q),
\intertext{and more generally, we have the cup product map}
\cup^n : \wedge^n [2^2] &\longrightarrow H^n(\Kg;\Q).
\end{align*}
Note that $\wedge^n [2^2]$ and $\Ker \cup^n$ are 
also $\Symp{Q}$-vector spaces for each $n \ge 2$. 
In this paper, we study the case of $n=2$. 
We put $\cup := \cup^2$, for simplicity. 
We will show the following.
\begin{llem}
For $g \ge 4$, the irreducible decomposition 
of $\wedge^2 [2^2]$ is given by
\[\wedge^2 [2^2] =
[431]+[42]+[3 2^2 1]+[321]+[31^3]+
[31]+[2^3]+[21^2]+[2]. \]
\end{llem}
\begin{tthm}
For $g \ge 4$, the kernel of the cup product map 
$\cup : \wedge^2 [2^2] \to H^2(\Kg;\Q)$ is 
\[[42]+[31^3]+[31]+[2^3]+[2], \]
which is, as an $\Symp{Q}$-vector space, 
isomorphic to the rational image of the fourth Johnson homomorphism 
$\tau_g (4)$.
\end{tthm}
\noindent
Here we write $+$ for the direct sum. Note that 
our theorem can be stated that we determine the kernel of 
the map induced on the second rational cohomology by $\tau_g (2)$. 
We will also treat the cases of $g=2,3$ in Section \ref{sec:g23}.

Before proving Theorem \ref{mainSg}, which corresponds to 
the case of a closed surface, 
we show a similar result 
for the case of a surface with a boundary in Section \ref{sec:compSg1}, 
since it is easier to handle by a technical reason. 
In each case, our task is divided into the following 
two parts. 
First, we will find some summands in Lemma \ref{irrSg1} (or 
that corresponding to each case) which belong to the kernel by 
using Stallings' exact sequence 
in \cite{st} together with Morita's description 
\cite{mo10, mo4} of Johnson's homomorphisms 
as a Lie algebra homomorphism. Then we show that the other summands 
actually survive in the second cohomology 
by constructing explicit cycles which 
come from abelian subgroups of the Johnson kernel and 
give non-trivial values by the Kronecker product. 

Note that the method we have explained now 
originated with Hain \cite{hain}, 
where he determined the kernel of the map induced 
on the second rational cohomology by the first Johnson 
homomorphism for the Torelli group. Our previous paper 
\cite{sa} treated the third cohomology of 
the Torelli group. 
Brendle-Farb \cite{bf} studied the second cohomology of 
the Torelli group and the Johnson kernel by using 
the Birman-Craggs-Johnson homomorphism, 
and 
Pettet \cite{pet} studied the second cohomology of the 
(outer-)automorphism group of a free group by using its first 
Johnson homomorphism. 

\vspace{10pt}

\section{Preliminaries}\label{sec:preliminaries}
\subsection{Surfaces and their mapping class groups}\label{subsec:MCG}
Let $H_1 (\Sg)$ be the first integral homology group of a closed oriented 
surface $\Sg$ of genus $g \ge 2$. 
$H_1 (\Sg)$ has a natural intersection form 
$\mu : H_1 (\Sg) \otimes H_1 (\Sg) \to \Z$ which 
is non-degenerate and skew symmetric. 
We fix a symplectic basis 
$\langle a_1,\ldots,a_g,b_1,\ldots,b_g \rangle$ of $H_1 (\Sg)$ 
with respect to $\mu$, namely 
\[\mu(a_i,a_j)=0,\quad \mu(b_i,b_j)=0,\quad
\mu(a_i,b_j)=\delta_{ij}.\]
The Poincar\'e duality gives a canonical isomorphism of 
$H_1 (\Sg)$ with its 
dual $H_1 (\Sg)^\ast =H^1 (\Sg)$, the 
first integral cohomology group of $\Sg$. In 
this isomorphism, $a_i$ (resp.\ $b_i$) $\in H_1 (\Sg)$ 
corresponds to $-b_i^\ast$ 
(resp.\ $a_i^\ast$) $\in H^1 (\Sg)$ 
where $\langle a_1^\ast ,\ldots,a_g^\ast ,b_1^\ast ,\ldots,b_g^\ast 
\rangle$ is the dual basis of $H^1 (\Sg)$. 
We use the same symbol $H$ for these canonically 
isomorphic abelian groups. 

We also use a compact oriented surface $\Sgb$ of genus $g$ 
with a connected boundary. $H_1 (\Sgb)$ can be naturally identified with 
$H$. The fundamental group $\pi_1 \Sgb$ of 
$\Sgb$, where we take a base point of $\Sgb$ on $\partial \Sgb$, 
is known to be a free group of rank $2g$. We write 
$\zeta \in \pi_1 \Sgb$ for the boundary loop of $\Sgb$. 
Then the fundamental group $\pi_1 \Sg$ of $\Sg$ is given by 
$\pi_1 \Sgb/\langle \zeta \rangle$ where $\langle \zeta \rangle$ is 
the normal closure of the subgroup generated by $\zeta$. 

Let $\Mg, \Mgp, \Mgb$ be the mapping class group of $\Sg$, 
of $\Sg$ relative to the base point, of $\Sgb$, respectively. 
They are related by the following exact sequences 
\[\begin{CD}
0 @>>> \Z @>>> \Mgb @>>> \Mgp @>>> 1, \\
1 @>>> \pi_1 \Sg @>>> \Mgp @>>> \Mg @>>> 1,
\end{CD}\]
where $\Z$ corresponds to the Dehn twist along a loop which is 
parallel to $\partial \Sgb$, and 
$\pi_1 \Sg$ is embedded in $\Mgp$ as spin-maps 
(see \cite[Theorem 4.3]{bi}).
The former sequence is a central extension. 

The natural action of $\Mg$ on $H$ gives the classical representation
\[\Mg \longrightarrow \Symp{Z},\]
and we also have similar ones for $\Mgp$ and $\Mgb$. The 
kernels of these representations are denoted by $\Ig$, $\Igp$ and 
$\Igb$, respectively and called the {\it Torelli group} for each case. 
Note that among $\Ig$, $\Igp$ and $\Igb$, 
we also have exact sequences similar to 
the above. Indeed the Dehn twist along $\partial \Sgb$ and spin-maps 
act on $H$ trivially. 

Let $\Kg$ (resp.\ $\Kgb$) be the subgroup of $\Mg$ (resp.\ $\Mgb$) generated 
by Dehn twists along bounding simple closed curves 
on $\Sg$ (resp.\ $\Sgb$). We define $\Kgp \subset \Mgp$ 
to be the image of $\Kgb$ by 
the map $\Mgb \to \Mgp$. Then we have 
\[\begin{CD}
0 @>>> \Z @>>> \Kgb @>>> \Kgp @>>> 1, \\
1 @>>> [\pi_1 \Sg ,\pi_1 \Sg] @>>> \Kgp @>>> \Kg @>>> 1,
\end{CD}\]
where the former sequence is the pull-back of 
the central extension of $\Mgp$, 
and the latter one follows from a result of Asada-Kaneko \cite{ak}.

\subsection{Johnson's homomorphisms}\label{subsec:Johnson}
In this subsection, 
we recall what we call Johnson's homomorphisms defined by 
Johnson \cite{jo, jo3} 
and Morita \cite{mo9, mo10, mo4, mo6}. 

By results of Dehn, Nielsen and many people, 
we have natural isomorphisms
\[\Mg \cong \, \Out_+ \pi_1 \Sg, \qquad
\Mgp \cong \, \Aut_+ \pi_1 \Sg, \qquad
\Mgb \cong \, \{ \varphi \in \Aut \pi_1 \Sgb \mid 
\varphi (\zeta)=\zeta \}, \]
where $\Out_+ \pi_1 \Sg := \Ker (\Out \pi_1 \Sg \to \Aut H_2 (\pi_1 \Sg))$, 
and $\Aut_+ \pi_1 \Sg$ is similar. 

For a group $G$, let $\{ \Gamma^k G \}_{k \ge 1}$ be 
the lower central series of $G$ 
inductively defined by $\Gamma^1 G=G$ and 
$\Gamma^{i} G =[\Gamma^{i-1} G ,G]$ for $i \ge 2$. 
By a general theory, \,$\{ (\Gamma^k G) /(\Gamma^{k+1} G)\}_{k \ge 1}$ 
forms a graded Lie algebra whose bracket map is induced from taking 
commutators. It is well known that the Lie algebra 
$\{ (\Gamma^k \pi_1 \Sgb) / (\Gamma^{k+1} \pi_1 \Sgb)\}_{k \ge 1}$ is 
isomorphic to the free Lie algebra $\Lgb= \{ \Lgb (k) \}_{k \ge 1}$ 
generated by $H$. 
Furthermore, by a result of Labute \cite{lab}, the Lie algebra 
$\{ (\Gamma^k \pi_1 \Sg) / (\Gamma^{k+1} \pi_1 \Sg)\}_{k \ge 1}$ is 
given by $\Lg := \Lgb/I$ where $I$ is the ideal of $\Lgb$ generated by 
$\omega_0:=\sum_{i=1}^g [a_i,b_i]$. 

The isomorphisms mentioned above induce the homomorphisms 
\begin{align*}
\sigma_{k} \ &: \hspace{3pt} \Mg \ \longrightarrow \Out 
\bigl( \pi_1 \Sg/(\Gamma^k \pi_1 \Sg) \bigr), \\
\sigma_{k,\ast} &: \Mgp \longrightarrow \Aut 
\bigl( \pi_1 \Sg/(\Gamma^k \pi_1 \Sg) \bigr), \\
\sigma_{k,1} &: \Mgb \longrightarrow \Aut 
\bigl( \pi_1 \Sgb/(\Gamma^k \pi_1 \Sgb) \bigr) ,
\end{align*}
for each $k \ge 2$, and we define filtrations of $\Mg, \Mgp, \Mgb$ by 
\begin{align*}
\Mg [1]:=& \, \Mg, \qquad \hspace{13pt} 
\Mg [k]:= \Ker \sigma_{k} \quad (k \ge 2),\\
\Mgp [1]:=& \, \Mgp, \qquad \Mgp [k]:= 
\Ker \sigma_{k,\ast} \quad (k \ge 2),\\
\Mgb [1]:=& \, \Mgb, \qquad \Mgb [k]:= 
\Ker \sigma_{k,1} \quad (k \ge 2). 
\end{align*}

For each $\varphi \in \Mgb [k+1]$ and $\gamma \in \pi_1 \Sgb$, 
we have $\varphi (\gamma) \gamma^{-1} \in \Gamma^{k+1} \pi_1 \Sgb$. 
This induces a map $\tau_{g,1} (k): \Mgb [k+1] \to 
\Hom (\pi_1 \Sgb, (\Gamma^{k+1} \pi_1 \Sgb)/(\Gamma^{k+2} \pi_1 \Sgb))= 
\Hom (H, \Lgb (k+1))$, 
and it is in fact a homomorphism. 

In \cite{mo10, mo4}, Morita showed the following. 
By taking commutators, we can endow 
$\{ \Mgb [k+1] /\Mgb [k+2] \}_{k \ge 1} = 
\{ \Im \tau_{g,1} (k) \}_{k \ge 1} =: \Im \tau_{g,1}$ 
with a Lie algebra structure. 
We can also endow 
$\Hom (H, \Lgb):=\{ \Hom (H, \Lgb (k+1)) \}_{k \ge 1}$ with 
a Lie algebra structure, so that 
$\tau_{g,1}:=\{ \tau_{g,1} (k) \}_{k \ge 1}$ becomes 
a Lie algebra inclusion of $\Im \tau_{g,1}$ into $\Hom (H, \Lgb)$. 
Moreover, $\Im \tau_{g,1}$ is contained in the Lie subalgebra 
$\DHgb = \{ \DHgb (k) \}_{k \ge 1}$ defined by 
\[\DHgb (k) := \Ker \left( \Hom (H, \Lgb (k+1)) \xrightarrow{\cong} 
H \otimes \Lgb (k+1) \xrightarrow{[\cdot,\cdot]} \Lgb (k+2) \right),\]
where the maps in the right hand side are 
given by the Poincar\'e duality and the bracket operation. 
A similar argument gives 
a homomorphism $\tau_{g,\ast} (k): \Mgp [k+1] \to \DHgp (k) \subset 
\Hom (H, \Lg (k+1))$, where
\[\DHgp (k) := \Ker \left( \Hom (H, \Lg (k+1)) \xrightarrow{\cong} 
H \otimes \Lg (k+1) \xrightarrow{[\cdot,\cdot]} \Lg (k+2) \right),\]
and the corresponding Lie algebra inclusion. 
Furthermore, Asada-Kaneko \cite{ak} showed that 
$\pi_1 \Sg \cap \Mgp [k+1] = \Gamma^{k} \pi_1 \Sg$. Hence we have 
an inclusion 
\[\Psi_k:(\Gamma^k \pi_1 \Sg)/(\Gamma^{k+1} \pi_1 \Sg) 
\cong \Lg (k) \hookrightarrow \DHgp (k).\]
If we set $\DHg (k) := \DHgp (k)/ \Lg (k)$, we obtain a homomorphism 
$\tau_g (k): \Mg [k+1] \to \DHg (k)$ (and the corresponding 
Lie algebra inclusion). 
We call the homomorphisms $\tau_g (k)$, $\tau_{g,\ast} (k)$, 
$\tau_{g,1} (k)$ the $k$-th Johnson homomorphism for each case. 
Note that $\tau_g (k)$ is $\Mg$-equivariant, where $\Mg$ acts on 
$\Mg [k+1]$ by conjugation 
and acts on the target through the classical representation 
$\Mg \to \Symp{Z}$. Similar results hold for $\tau_{g,\ast} (k)$ and 
$\tau_{g,1} (k)$. 

We have $\Mgb[2]=\Igb$, $\Mgp[2]=\Igp$ and $\Mg[2]=\Ig$ by definition. 
Johnson \cite{jo1} showed that $\Mgb[3]=\Kgb$ and $\Mg[3]=\Kg$. 
Combining the fact that the first Johnson homomorphisms for 
$\Igb$ and $\Igp$ have the same target $\wedge^3 H$, we can see 
that $\Mgp[3]=\Kgp$.

\subsection{The representation theory of $\Symp{Q}$}
Here we summarize the notation and 
general facts concerning 
the representation theory of $\Symp{Q}$ 
from \cite{fh}, \cite{hain} and \cite{mo6}. 
First we consider the Lie group $\Symp{C}$ and 
its Lie algebra $\symp{C}$. By a general 
theory of the representation, finite dimensional 
representations of $\Symp{C}$ coincide with those of $\symp{C}$, 
and their common irreducible representations (up to isomorphisms) are 
parameterized by Young diagrams whose numbers of 
rows are less than or equal to $g$. These representations 
are all defined over $\Q$ so that we can consider them as irreducible 
representations of $\Symp{Q}$ and $\symp{Q}$. 
We follow the notation in \cite{mo6} to describe Young diagrams. 
For example, the trivial representation $\Q$ is denoted by $[0]$ and 
the fundamental representation $\Hq:= H \otimes \Q$ is denoted by $[1]$. 
We fix a symplectic basis $\langle a_1,\ldots,a_g,b_1,\ldots,b_g \rangle$ 
of $\Hq$ with respect to the non-degenerate skew symmetric bilinear form, 
denoted by $\mu$ again, on $\Hq$ induced from 
the intersection form $\mu$ on $H$. 
In general, the Young diagram $[n_1 n_2 \cdots n_l]$, where $n_i$ are 
integers satisfying $n_1 \ge n_2 \ge \cdots \ge n_l \ge 1$ 
and $l \le g$, corresponds to the $\Symp{Q}$-vector space $V$ 
given as follows. Let $[m_1 m_2 \cdots m_k]$ be the Young diagram 
obtained by transposing $[n_1 n_2 \cdots n_l]$. Then $V$ is explicitly 
defined to be the irreducible $\Symp{Q}$-subspace of 
\[(\wedge^{m_1} \Hq) \otimes (\wedge^{m_2} \Hq) \otimes 
\cdots \otimes (\wedge^{m_k}\Hq) \]
containing the vector
\[(a_1 \wedge a_2 \wedge \cdots \wedge a_{m_1}) \otimes 
(a_1 \wedge a_2 \wedge \cdots \wedge a_{m_2}) \otimes 
\cdots \otimes 
(a_1 \wedge a_2 \wedge \cdots \wedge a_{m_k}), 
\]
which is called 
the {\it highest weight vector of $[n_1 n_2 \cdots n_l]$}. 

We define the following elements 
$X_{i,j}, Y_{i,j}$ $(i \neq j)$ and $U_i, V_i$ of $\symp{Q}$ 
characterized by their actions on $\Hq$ as follows:
\[\begin{array}{rclrcl}
X_{i,j}(a_k) &=& \delta_{jk}a_i, & X_{i,j}(b_k) &=& -\delta_{ik}b_j, \\
Y_{i,j}(a_k) &=& 0, & Y_{i,j}(b_k) &=& \delta_{ik}a_j + \delta_{jk}a_i, \\
U_{i}(a_k) &=& 0, & U_{i}(b_k) &=& \delta_{ik}a_i, \\
V_{i}(a_k) &=& \delta_{ik}b_i, & V_{i}(b_k) &=& 0.
\end{array}\]

We will frequently make use of the following $\Symp{Q}$-(and $\symp{Q}$-)
equivariant homomorphisms. 
\begin{itemize}
\item[1)] \; The \textit{contraction} 
$C_k^{(i,j)} :\otimes^k \Hq \rightarrow \otimes^{k-2} \Hq$ 
$(1 \le i < j \le k)$ is given by
\[C_k^{(i,j)} (x_1 \otimes \cdots \otimes x_k)=
\mu(x_i,x_j) x_1 \otimes \cdots \otimes \widehat{x_i} 
\otimes \cdots \otimes \widehat{x_j} \otimes \cdots \otimes x_k,\]
where $\widehat{x_i}$ means excluding $x_i$, and 
$\otimes^0 \Hq=\Q$ is the trivial representation.

\item[2)] \;
For each symbol 
\[\sigma=(i(1,1),\ldots,i(1,k_1))(i(2,1),\ldots,i(2,k_2))\cdots
(i(l,1),\ldots,i(l,k_l)),\]
where $k=\sum_{j=1}^l k_j$ and $\{ i(1,1),\ldots,
i(1,k_1),i(2,1),\ldots,i(l,k_l)\}=\{1,2,\ldots,k \}$, the projection 
$p_k^\sigma : \otimes^k \Hq \rightarrow 
(\wedge^{k_1} \Hq) \otimes (\wedge^{k_2} \Hq) \otimes \cdots \otimes 
(\wedge^{k_l} \Hq)$ is given by 
\[p_k^\sigma (x_1 \otimes \cdots \otimes x_k)=(x_{i(1,1)} \wedge 
\cdots \wedge x_{i(1,k_1)}) \otimes \cdots \otimes (x_{i(l,1)} \wedge 
\cdots \wedge x_{i(l,k_l)}).\]
For example, $p_3^{(1,3)(2)}(x_1 \otimes x_2 \otimes x_3)=
(x_1 \wedge x_3) \otimes x_2$. 

\item[3)] \; The canonical inclusion 
$\iota_k : \Lgb^\Q (k) := \Lgb (k) \otimes \Q 
\hookrightarrow \Hq^{\otimes k}$, 
where $\Lgb (k)$ is the degree $k$ part 
of the free Lie algebra $\Lgb$ generated by $H$, 
is inductively defined by 
replacing the bracket $[X,Y]$ by $X \otimes Y - Y \otimes X$.

\end{itemize}

\subsection{Johnson's homomorphisms via 
the representation theory of $\Symp{Q}$}\label{subsec:Johnson-and-rep} 
Let $\tau_{g}^\Q (k)$ denote the Johnson homomorphism $\tau_{g} (k)$ 
tensored by $\Q$, namely 
\[\tau_{g}^\Q (k): (\Mg [k+1]/(\Gamma^2 \Mg [k+1])) \otimes \Q 
\longrightarrow \DHg^\Q (k):= \DHg (k) \otimes \Q.\]
As mentioned in Section \ref{subsec:Johnson}, $\tau_{g}^\Q (k)$ is 
$\Mg$-equivariant, and in particular, $\Im \tau_{g}^\Q (k)$ is 
an $Sp(2g,$ $\Z)$-vector space. Moreover, it turns out that 
$\Im \tau_{g}^\Q (k)$ is in fact an $\Symp{Q}$-vector space 
by Asada-Nakamura \cite[Lemma 2.2.8]{an}. 
Similar results hold for $\tau_{g,\ast}^\Q (k):=\tau_{g,\ast} (k)
\otimes \Q$ and $\tau_{g,1}^\Q (k):=\tau_{g,1} (k) \otimes \Q$. 
By results of Johnson \cite{jo} for $k=1$, 
Hain \cite{hain} and Morita \cite{mo9} for $k=2$, 
Hain \cite{hain} and
Asada-Nakamura \cite{an} for $k=3$, we have the following. 
\begin{align*}
\Im \tau_{g,1}^\Q (1) &= \Im \tau_{g,\ast}^\Q (1) = \DHgb^\Q (1) =
\DHgp^\Q (1) = [1^3]+[1] = \wedge^3 \Hq, \\
\Im \tau_{g}^\Q (1) \ &= \DHg^\Q (1) \ \hspace{1pt} = [1^3] = 
(\wedge^3 \Hq) / \Hq, \\
\Im \tau_{g,1}^\Q (2) &= \DHgb^\Q (2) = [2^2]+[1^2]+[0], \\
\Im \tau_{g,\ast}^\Q (2) &= \DHgp^\Q (2) = [2^2]+[1^2], \\
\Im \tau_{g}^\Q (2) \ &= \DHg^\Q (2) \ \hspace{1pt} = [2^2] \\
\Im \tau_{g,1}^\Q (3) &= \Im \tau_{g,\ast}^\Q (3) = [31^2]+[21] 
\subset \DHgb^\Q (3) = \DHgp^\Q (3) = [31^2]+[21]+[3], \\
\Im \tau_{g}^\Q (3) \ &= [31^2] \subset \DHg^\Q (3) = [31^2]+[3] 
\intertext{for $g \ge 3$. 
Moreover, in \cite{mo6}, Morita announced that}
\Im \tau_{g,1}^\Q (4) &= \Im \tau_{g,\ast}^\Q (4) = 
[42]+[31^3]+2[31]+[2^3]+[21^2]+2[2], \\
\Im \tau_{g}^\Q (4) \ &= [42]+[31^3]+[31]+[2^3]+[2]
\intertext{in}
\DHgb^\Q (4) &= [42]+[31^3]+2[31]+[2^3]+2[21^2]+3[2], \\
\DHgp^\Q (4) &= [42]+[31^3]+2[31]+[2^3]+2[21^2]+2[2], \\
\DHg^\Q (4) \ &= [42]+[31^3]+[31]+[2^3]+[21^2]+[2]
\end{align*}
for $g \ge 4$. 
\begin{remark}
Hain showed in \cite{hain} that 
as Lie algebras, $\Im \tau_{g,1}^\Q, \Im \tau_{g,\ast}^\Q, 
\Im \tau_{g}^\Q$ are generated by their degree 1 parts, and 
that $\Im \tau_{g,1}^\Q (k) \cong \Im \tau_{g,\ast}^\Q (k)$ for $k \ge 3$. 
\end{remark}
\begin{remark}
A part of our proof of Theorems \ref{mainSg1}, \ref{mainSg} in 
later sections has overlaps with Morita's proof \cite{mop} of 
the determinations of $\Im \tau_{g,1}^\Q (4) \cong 
\Im \tau_{g,\ast}^\Q (4)$ and $\Im \tau_{g}^\Q (4)$. 
More precisely, we will prove that the summands displayed above 
are contained in $\Im \tau_{g,1}^\Q (4)$, 
$\Im \tau_{g,\ast}^\Q (4)$ and $\Im \tau_{g}^\Q (4)$ by using 
brackets of two elements in $\Im \tau_{g,1}^\Q (2)$, 
$\Im \tau_{g,\ast}^\Q (2)$ and $\Im \tau_{g}^\Q (2)$. As a result, 
we can see that the degree 4 parts of the rational images of Johnson's 
homomorphisms are generated by their degree 2 parts. 
\end{remark}

The bracket operation of $\DHgb$ is explicitly 
given in \cite{mo10, mo4}. 
In this paper, however, we use an alternative description given by 
Garoufalidis-Levine \cite{gl}, Levine \cite{le3, le4}, 
which will be easier to handle. 
We recall the following Lie algebra of labeled unitrivalent trees. 
Let $\At{k}$ be the abelian group generated by unitrivalent trees 
with $k+2$ univalent vertices labeled by elements of $H$ and 
a cyclic order of each trivalent vertex modulo relations of AS and 
IHX together with linearity of labels. We can endow $\At{} := 
\{ \At{k} \}_{k \ge 1}$ with a bracket 
operation 
\[ [\ \cdot\ ,\ \cdot\ ]:\At{k} \otimes \At{l} \longrightarrow \At{k+l}\]
given below, and $\At{}$ becomes a quasi Lie algebra. 
For labeled trees $T_1, T_2 \in \At{}$, we define 
\[[T_1,T_2]:= \sum_{i,j} \mu (a_i,b_j) T_1 \ast_{i,j} T_2\]
where the sum is taken over all pairs of a univalent vertex of $T_1$, 
labeled by $a_i$, and one of $T_2$, labeled by $b_j$, and 
$T_1 \ast_{i,j} T_2$ is the tree given by welding $T_1$ and $T_2$ at 
the pair. We define a map $\eta_k: \At{k} \to H \otimes \Lgb (k+1)$ by 
\[\eta_k (T) := \sum_{v} a_v \otimes T_v,\]
where the sum is over all univalent vertices of $T$, and for each univalent 
vertex $v$, $a_v$ denotes the label of $v$ and $T_v$ denotes 
the rooted labeled planar binary tree obtained by removing 
the label $a_v$ and considering $v$ to be an unlabeled root, 
which can be regarded as 
an element of $\Lgb (k+1)$ by a standard method. It is shown 
that $\eta:=\{\eta_k\}_{k \ge 1}:\At{} \to H \otimes \Lgb$ 
is a quasi Lie algebra homomorphism 
and $\Im \eta \subset \DHgb$. Moreover $\eta \otimes \Q : 
\At{} \otimes \Q \to \DHgb^\Q$ becomes an isomorphism of Lie 
algebras. For more details, 
see \cite{gl}, \cite{le3, le4} and their references. 
In what follows, we identify $\At{} \otimes \Q$ with 
$\DHgb^\Q$ by $\eta \otimes \Q$. 

Using $\At{} \otimes \Q$, 
we now give a graphical description of the map 
\[\Psi_k:\Lg^\Q (k):= \Lg (k) \otimes \Q \hookrightarrow 
\Im \tau_{g,\ast}^\Q (k) \subset \DHgp^\Q (k),\] 
which was mentioned in Section \ref{subsec:Johnson} 
and is explicitly given by 
\[\Lg^\Q (k) \ni X \, \mapsto \, 
\sum_{i=1}^g \big( a_i \otimes [b_i,X] - b_i \otimes [a_i,X] \big) 
\in H \otimes \Lg^\Q (k+1),\]
as follows. For each rooted labeled planar binary tree $T$ as an 
element of $\Lgb^\Q (k)$, we can construct an element of 
$\At{k} \otimes \Q \cong \DHgb^\Q (k)$ by gluing $T$ to the rooted 
labeled planar binary tree $T_{\omega_0} \in \Lgb^\Q (2) \cong 
\wedge^2 \Hq$ corresponding to $\omega_0 =\sum_{i=1}^g a_i \wedge b_i$ 
at their roots as depicted in Figure 1. 

\vspace*{20pt}
%\hspace*{-50pt}
%\input{graph1.tex}
%WinTpicVersion3.08
\unitlength 0.1in
\begin{picture}( 58.2000, 12.9500)(  1.8000,-15.1500)
% CIRCLE 1 0 1 0
% 4 805 805 1005 805 1005 805 1205 805
% 
\special{pn 13}%
\special{sh 0.300}%
\special{ar 806 806 200 200  0.0000000 6.2831853}%
% STR 2 0 3 0
% 3 805 705 805 805 5 0
% $X$
\put(8.0500,-8.0500){\makebox(0,0){$X$}}%
% STR 2 0 3 0
% 3 405 705 405 805 5 0
% $T =$
\put(4.0500,-8.0500){\makebox(0,0){$T =$}}%
% LINE 1 0 3 0
% 2 1005 805 1405 805
% 
\special{pn 13}%
\special{pa 1006 806}%
\special{pa 1406 806}%
\special{fp}%
% STR 2 0 3 0
% 3 2210 705 2210 805 5 0
% $T_{\omega_0} = \displaystyle\sum_{i=1}^g$
\put(22.1000,-8.0500){\makebox(0,0){$T_{\omega_0} = 
\displaystyle\sum_{i=1}^g$}}%
% LINE 1 0 3 0
% 6 2710 805 3110 805 3310 405 3110 805 3310 1205 3110 805
% 
\special{pn 13}%
\special{pa 2710 806}%
\special{pa 3110 806}%
\special{fp}%
\special{pa 3310 406}%
\special{pa 3110 806}%
\special{fp}%
\special{pa 3310 1206}%
\special{pa 3110 806}%
\special{fp}%
% STR 2 0 3 0
% 3 3410 205 3410 305 5 0
% $b_i$
\put(34.1000,-3.0500){\makebox(0,0){$b_i$}}%
% STR 2 0 3 0
% 3 3410 1205 3410 1305 5 0
% $a_i$
\put(34.1000,-13.0500){\makebox(0,0){$a_i$}}%
% CIRCLE 1 0 1 0
% 4 5200 800 5400 800 5400 800 5600 800
% 
\special{pn 13}%
\special{sh 0.300}%
\special{ar 5200 800 200 200  0.0000000 6.2831853}%
% STR 2 0 3 0
% 3 4500 710 4500 810 5 0
% $\Phi_k (T)=\displaystyle\sum_{i=1}^g$
\put(45.0000,-8.1000){\makebox(0,0){$\Phi_k (T)=
\displaystyle\sum_{i=1}^g$}}%
% LINE 1 0 3 0
% 6 5400 805 5800 805 6000 405 5800 805 6000 1205 5800 805
% 
\special{pn 13}%
\special{pa 5400 806}%
\special{pa 5800 806}%
\special{fp}%
\special{pa 6000 406}%
\special{pa 5800 806}%
\special{fp}%
\special{pa 6000 1206}%
\special{pa 5800 806}%
\special{fp}%
% STR 2 0 3 0
% 3 6100 205 6100 305 5 0
% $b_i$
\put(61.0000,-3.0500){\makebox(0,0){$b_i$}}%
% STR 2 0 3 0
% 3 6100 1205 6100 1305 5 0
% $a_i$
\put(61.0000,-13.0500){\makebox(0,0){$a_i$}}%
% STR 2 0 3 0
% 3 5200 700 5200 800 5 0
% $X$
\put(52.0000,-8.0000){\makebox(0,0){$X$}}%
% CIRCLE 2 0 2 0
% 4 1405 805 1455 805 1455 805 1655 805
% 
\special{pn 8}%
\special{sh 0}%
\special{ar 1406 806 50 50  0.0000000 6.2831853}%
% CIRCLE 2 0 2 0
% 4 2710 805 2760 805 2760 805 2960 805
% 
\special{pn 8}%
\special{sh 0}%
\special{ar 2710 806 50 50  0.0000000 6.2831853}%
% STR 2 0 3 0
% 3 3850 710 3850 810 5 0
% $\Longrightarrow$
\put(38.5000,-8.1000){\makebox(0,0){$\Longrightarrow$}}%
% STR 2 0 3 0
% 3 3100 1500 3100 1600 5 0
% The map $\Phi_k:\Lgb^\Q (k) \to \DHgb^\Q (k)$
\put(31.0000,-16.0000){\makebox(0,0){{\footnotesize Figure 1. The map 
$\Phi_k:\Lgb^\Q (k) \to \DHgb^\Q (k)$}}}%
% STR 2 0 3 0
% 3 1400 830 1400 930 5 0
% {\tiny root}
\put(14.0000,-9.3000){\makebox(0,0){{\tiny root}}}%
% STR 2 0 3 0
% 3 2700 830 2700 930 5 0
% {\tiny root}
\put(27.0000,-9.3000){\makebox(0,0){{\tiny root}}}%
\end{picture}%

\vspace{15pt}

It can be easily checked that this construction gives an 
$\Symp{Q}$-equivariant homomorphism 
$\Phi_k:\Lgb^\Q (k) \to \DHgb^\Q (k)$, and 
it induces the desired map 
$\Psi_k:\Lg^\Q (k) \hookrightarrow \DHgp^\Q (k)$ by using 
Labute's result \cite{lab}. 

\vspace{10pt}

\section{Description of the second Johnson homomorphism}
\label{sec:secondjohnson}
In this paper, we mainly concern the second Johnson homomorphism and 
the map induced on the rational cohomology. 
Now we see the second Johnson homomorphism more closely. 
Recall that 
\[\begin{CD}
\Kgb @>\tau_{g,1}^\Q(2)>> \DHgb^\Q (2) @= [2^2]+[1^2]+[0] \\
@VVV @VVV @VV\mathrm{proj}V \\
\Kgp @>\tau_{g,\ast}^\Q(2)>> \DHgp^\Q (2) @= [2^2]+[1^2] \\
@VVV @VVV @VV\mathrm{proj}V \\
\Kg @>\tau_{g}^\Q(2)>> \DHg^\Q (2) @= [2^2] 
\end{CD}\]
for $g \ge 2$. Note that by passing to the duals, 
we have $\Symp{Q}$-equivariant inclusions 
\[[2^2] \hookrightarrow [2^2]+[1^2] \hookrightarrow [2^2]+[1^2]+[0]\]
which are unique up to scalars. 

The summands $[0]$ and $[1^2]$ in $\DHgb^\Q (2)$ can be detected by 
the composite of maps 
\begin{align*}
q_{[1^2]} :& \ \DHgb^\Q (2) \hookrightarrow H \otimes 
\Lgb^\Q (3) \stackrel{1 \otimes \iota_3}{\hookrightarrow} 
\Hq^{\otimes 4} \xrightarrow{C_4^{(1,2)}} \Hq^{\otimes 2} 
\xrightarrow{p_2^{(1,2)}} \wedge^2 \Hq, \\
q_{[0]} :& \ \DHgb^\Q (2) \hookrightarrow 
H \otimes \Lgb^\Q (3) \stackrel{1 \otimes \iota_3}{\hookrightarrow} 
\Hq^{\otimes 4} \xrightarrow{C_2^{(1,2)} \circ C_4^{(1,2)}} \Q=[0],
\end{align*}
\noindent
respectively. On the other hand, by using the map 
$\Phi_2:\Lgb^\Q (2)=\wedge^2 \Hq =[1^2]+[0] 
\to \DHgb^\Q (2)$ mentioned in Section 
\ref{subsec:Johnson-and-rep}, we can obtain the highest 
weight vectors $v_{[2^2]}$, $v_{[1^2]}$ and $v_{[0]}$ 
of $[2^2]$, $[1^2]$ and $[0]$ defined by 
\[v_{[2^2]} = T^H(a_1,a_2,a_1,a_2), \quad 
v_{[1^2]} = \sum_{i=1}^g \ T^H(a_1,a_2,a_i,b_i), \quad
v_{[0]} = \sum_{i,j=1}^g \ T^H(a_i,b_i,a_j,b_j),\]
where we put

\vspace*{15pt}
\hspace*{50pt}
%WinTpicVersion3.08
\unitlength 0.1in
\begin{picture}( 26.5500, 10.0000)( -0.5500,-12.1500)
% LINE 1 0 3 0
% 10 1800 800 1600 400 1600 1200 1800 800 1800 800 
%2400 800 2400 800 2600 400 2600 1200 2400 800
% 
\special{pn 13}%
\special{pa 1800 800}%
\special{pa 1600 400}%
\special{fp}%
\special{pa 1600 1200}%
\special{pa 1800 800}%
\special{fp}%
\special{pa 1800 800}%
\special{pa 2400 800}%
\special{fp}%
\special{pa 2400 800}%
\special{pa 2600 400}%
\special{fp}%
\special{pa 2600 1200}%
\special{pa 2400 800}%
\special{fp}%
% STR 2 0 3 0
% 3 1500 200 1500 300 5 0
% $a$
\put(15.0000,-3.0000){\makebox(0,0){$a$}}%
% STR 2 0 3 0
% 3 1500 1200 1500 1300 5 0
% $b$
\put(15.0000,-13.0000){\makebox(0,0){$b$}}%
% STR 2 0 3 0
% 3 2700 1200 2700 1300 5 0
% $c$
\put(27.0000,-13.0000){\makebox(0,0){$c$}}%
% STR 2 0 3 0
% 3 2700 200 2700 300 5 0
% $d$
\put(27.0000,-3.0000){\makebox(0,0){$d$}}%
% STR 2 0 3 0
% 3 800 700 800 800 5 0
% $T^H(a,b,c,d) \ :=$
\put(8.0000,-8.0000){\makebox(0,0){$T^H(a,b,c,d) \ :=$}}%
% STR 2 0 3 0
% 3 3300 700 3300 800 5 0
% $\in \ \DHgb^\Q (2).$
\put(33.0000,-8.0000){\makebox(0,0){$\in \ \DHgb^\Q (2).$}}%
\end{picture}%

\vspace{15pt}

In \cite{mo9}, Morita computed the value of the Dehn twist 
$\psi_h \in \Kgb$ 
along the simple closed curve $\gamma_h$ bounding the surface 
$\Sigma_{h,1}$ of the standard position by $\tau_{g,1}^\Q (2)$, and 
he obtained that 
\begin{align*}
\tau_{g,1}^\Q (2) (\psi_h) &= \sum_{i,j=1}^h \big\{ 
a_i \otimes [[a_j,b_j],b_i]-b_i \otimes [[a_j,b_j],a_i] \big\} \\
&= -\frac{1}{2} \ \sum_{i,j=1}^h \ T^H (a_i,b_i,a_j,b_j).
\end{align*}
In particular, we have $\tau_{g,1}^\Q (2) (\psi_1^{-2})=
T^H (a_1,b_1,a_1,b_1)$.

Let $\overline{\psi_1^{-2}} \in \Kg$ be the image of 
$\psi_1^{-2} \in \Kgb$ under 
the projection $\Kgb \to \Kg$. To compute $\tau_{g}^\Q (2) 
\bigl( \overline{\psi_1^{-2}} \bigr)$, we need to project the vector 
$\tau_{g,1}^\Q (2) (\psi_1^{-2})=T^H(a_1,b_1,a_1,b_1)$ onto 
$[2^2] \subset [2^2]+[1^2]+[0]$. 
By direct computations, we can see that 
\begin{align*}
&q_{[1^2]}(T^H(a_1,b_1,a_1,b_1)) = 12 a_1 \wedge b_1, 
&&q_{[0]}(T^H(a_1,b_1,a_1,b_1)) = 12, \\
&q_{[1^2]}(\Phi_2 (a_1 \wedge b_1)) = (4g+4) a_1 \wedge b_1 +4 \omega_0, 
&&q_{[0]}(\Phi_2 (a_1 \wedge b_1)) = 8g+4,\\
&q_{[1^2]}(\Phi_2 (\omega_0)) = q_{[1^2]}(v_{[0]})  = (8g+4) \omega_0, 
&&q_{[0]}(\Phi_2 (\omega_0)) = q_{[0]}(v_{[0]}) = 8g^2+4g,
\end{align*}
and therefore we have
\begin{align*}
\tau_{g}^\Q (2) \bigl( \overline{\psi_1^{-2}} \bigr) 
&= T^H (a_1,b_1,a_1,b_1)-\frac{3}{g+1}\Phi_2 (a_1 \wedge b_1)
+\frac{3}{(2g+1)(g+1)} \Phi_2 (\omega_0) \\
&\in [2^2]=\DHg (2).
\end{align*}

Consequently, we have the following.
\begin{prop}\label{image}
$(1)$ {\rm(Morita \cite{mo9})} $\Im \tau_{g,1}^\Q (2)$ contains the vector 
$T^H (a_1,b_1,a_1,b_1)$. \\
$(2)$ \ $\Im \tau_{g}^\Q (2)$ contains the vector 
\[T^H (a_1,b_1,a_1,b_1)
-\frac{3}{g+1}\sum_{i=1}^g T^H (a_1,b_1,a_i,b_i)
+\frac{3}{(2g+1)(g+1)}\sum_{i,j=1}^g T^H(a_i,b_i,a_j,b_j).\]
\end{prop}

\vspace{10pt}

\section{Computation for the case of a surface with a boundary}
\label{sec:compSg1}
\subsection{Statement for $\Kgb$}\label{subsec:statementSg1}
We begin our proof of the main theorem. First, we 
treat the case of $\Kgb$, whose 
computational results will be used again 
for proofs of the cases of $\Kgp$ and $\Kg$. 

As seen in Section \ref{subsec:Johnson-and-rep}, 
we have $\Im \tau_{g,1}^\Q (2) = \DHgb^\Q (2) \cong [2^2]+[1^2]+[0]$, 
so that an injection $[2^2]+[1^2]+[0] \hookrightarrow H^1 (\Kgb;\Q)$ 
is obtained. We now consider the cup product map
\[\cup : \wedge^2 ([2^2]+[1^2]+[0]) \longrightarrow H^2 (\Kgb;\Q).\]
Note that 
\[\wedge^2 ([2^2]+[1^2]+[0]) 
\cong \wedge^2 [2^2] + [2^2] \otimes [1^2] + \wedge^2 [1^2] + 
[2^2] \otimes [0] + [1^2] \otimes [0].\]
\begin{lem}\label{irrSg1}
For $g \ge 4$, the irreducible decompositions 
of $\wedge^2 [2^2]$, $[2^2] \otimes [1^2]$, $\wedge^2 [1^2]$, 
$[2^2] \otimes [0]$ and $[1^2] \otimes [0]$ are 
given by the following table. 

\[\begin{array}{c}
\begin{array}{|c||c|c|c|c|c|}
\hline
 & \wedge^2 [2^2] & [2^2] \otimes [1^2]  
 & \wedge^2 [1^2] & [2^2] \otimes [0]  
 & [1^2] \otimes [0]  \\
\hline
 [431]  &  1  & & & & \\
\hline
 [42]  &  1  & & & & \\
\hline
 [32^2 1]  &  1  & & & & \\
\hline
 [3 2 1]  &  1  &  1  & & & \\
\hline
 [3 1^3]  &  1  & & & & \\
\hline
 [3 1]  &  1  &  1  & & & \\
\hline
 [2^3]  &  1  & & & & \\
\hline
 [2 1^2]  &  1  &  1  &  1  & & \\
\hline
 [2]  &  1  & &  1  & & \\
\hline
 [3^2]  & &  1  & & & \\
\hline
 [2^2 1^2]  & &  1  & & & \\
\hline
 [2^2]  & &  1  &  &  1  & \\
\hline
 [1^2]  & &  1  & & &  1  \\
\hline
\end{array}\vspace{2mm}\\
\mbox{\footnotesize {\rm Table 1}. Irreducible decompositions 
of $\Symp{Q}$-modules for $g \ge 4$} 
\end{array}\]

\noindent
In particular, we have
\begin{align*}
\wedge^2 ([2^2]+[1^2]+[0]) \cong &\ ([431]+[42]+[3 2^2 1]+[321]+[31^3]+
[31]+[2^3]+[21^2]+[2]) \\
&+ ([321]+[31]+[21^2]+[3^2]+[2^2 1^2]+[2^2]+[1^2]) \\
&+ ([21^2]+[2]) + [2^2] + [1^2].
\end{align*}
\end{lem}
\begin{proof}
By a general theory of the representation (see \cite[Remark 6.2]{hain}), 
it is observed that the irreducible decompositions we now consider 
become stable for $g \ge 4$. The computations in the next subsections 
show that the summands in the table are actually contained. To show that 
they are all, it suffices to check that the total dimension of the summands 
coincides with the dimension of $\wedge^2 \DHgb^\Q (2)$ in the case $g=4$ by 
using Weyl's character formula (see \cite[Section 24.2]{fh}). We now 
omit the details. 
\end{proof}
The first main theorem of this section is the following.
\begin{thm}\label{mainSg1}
For $g \ge 4$, the kernel of the cup product map 
$\cup : \wedge^2 ([2^2]+[1^2]+[0]) \to H^2(\Kgb;\Q)$ is 
\[[42]+[31^3]+2[31]+[2^3]+[21^2]+2[2], \]
which is, as an $\Symp{Q}$-vector space, 
isomorphic to $\Im \tau_{g,1}^\Q (4)$.
\end{thm}
\noindent
This theorem follows from Lemmas 
\ref{lem:lowerSg1}, \ref{lem:upperSg1} 
in the next two subsections .

\subsection{A lower bound of the kernel}
\label{subsec:lowerSg1}
Here we show that the summands in Theorem \ref{mainSg1} are 
contained in the kernel of the cup product map. 
Note that the map 
$\cup : \wedge^2 ([2^2]+[1^2]+[0]) \to H^2(\Kgb;\Q)$ 
is nothing other than the homomorphism 
\begin{align*}
\tau_{g,1}(2)^\ast &:H^2 (\Im \tau_{g,1} (2);\Q) \longrightarrow 
H^2 (\Kgb;\Q).
\intertext{Hence, by passing to the dual, 
our task is equivalent to observing the cokernel of the map}
\tau_{g,1}(2)_\ast &:H_2 (\Kgb;\Q) \longrightarrow 
H_2 (\Im \tau_{g,1} (2);\Q).
\end{align*}

By applying Stallings' exact sequence \cite{st} to the group 
extension 
\[1 \longrightarrow \Gamma^2 (\Kgb) \longrightarrow \Kgb 
\longrightarrow H_1 (\Kgb) \longrightarrow 1\]
and observing the homomorphisms, 
we obtain the exact sequence 
\[H_2 (\Kgb;\Q) \longrightarrow \wedge^2 H_1 (\Kgb;\Q) 
\xrightarrow{[\cdot,\cdot]} 
((\Gamma^2 \Kgb) / (\Gamma^{3} \Kgb)) \otimes \Q \longrightarrow 0\]
where the first map is the coproduct on the rational homology and 
the second one is the Lie bracket 
\[[\ \cdot\ ,\ \cdot\ ]:\wedge^2 H_1 (\Kgb) = 
\wedge^2 ( (\Gamma^1 \Kgb) / (\Gamma^{2} \Kgb) ) 
\longrightarrow 
(\Gamma^2 \Kgb) / (\Gamma^{3} \Kgb).\]

Since $\Gamma^2 \Kgb \subset \Mgb[5]$ and 
$\Gamma^3 \Kgb \subset \Mgb[7]$, we have a natural map 
\[i:(\Gamma^2 \Kgb) / (\Gamma^{3} \Kgb) \longrightarrow \Mgb[5]/\Mgb[6].\] 
By considering the fact that 
$\tau_{g,1}:\{ \Mgb [k+1] /\Mgb [k+2] \}_{k \ge 1} \to 
\Im \tau_{g,1} \subset \DHgb$ is a 
Lie algebra homomorphism as mentioned in Section \ref{subsec:Johnson}, 
we have the commutative diagram
\[\begin{CD}
\wedge^2 H_1 (\Kgb) @>{[\cdot,\cdot]}>> 
(\Gamma^2 \Kgb) / (\Gamma^{3} \Kgb) @>i>> \Mgb[5]/\Mgb[6]\\
@V{\wedge^2 \tau_{g,1} (2)}VV @. @V{\cong}V{\tau_{g,1} (4)}V \\
\wedge^2 \Im \tau_{g,1} (2) @>{[\cdot,\cdot]}>> 
\Im \tau_{g,1} (4) @= \Im \tau_{g,1} (4).
\end{CD}\]
Consequently, we obtain the commutative diagram
\[\begin{CD}
H_2 (\Kgb;\Q) @>>> \wedge^2 H_1 (\Kgb;\Q) @>[\cdot,\cdot]>>
((\Gamma^2 \Kgb) / (\Gamma^{3} \Kgb)) \otimes \Q \\
@V{\tau_{g,1} (2)_\ast}VV @V{\wedge^2 \tau_{g,1}^\Q (2)}VV 
@VV{\tau_{g,1}^\Q (4) \circ i}V \\
H_2 (\Im \tau_{g,1} (2);\Q) @>\cong>> \wedge^2 \Im \tau_{g,1}^\Q (2) 
@>[\cdot,\cdot]>> \Im \tau_{g,1}^\Q (4)
\end{CD}\]
whose upper row is exact. Hence we can find some summands in 
the cokernel of $\tau_{g,1}(2)_\ast :H_2 (\Kgb;\Q) \to 
H_2 (\Im \tau_{g,1} (2);\Q)$ by observing 
the image of 
\[[\ \cdot\ ,\ \cdot\ ] \ :\ \wedge^2 \Im \tau_{g,1}^\Q (2) 
\longrightarrow \Im \tau_{g,1}^\Q (4).\]
\begin{remark}\label{referee}
It is not known whether $H_1 (\Kgb;\Q)$ is finite dimensional or not. 
If it is finite dimensional, 
we could use Sullivan's exact sequence \cite{sullivan} in the 
above argument as in Pettet's one \cite[Section 3.2]{pet}.
\end{remark}

\newpage

\begin{lem}\label{lem:lowerSg1}
For $g \ge 4$, the image of the Lie bracket 
$[\cdot,\cdot] : \wedge^2 \Im \tau_{g,1}^\Q (2) \to 
\Im \tau_{g,1}^\Q (4)$ contains 
\[[42]+[31^3]+2[31]+[2^3]+[21^2]+2[2].\]
Hence, these summands are contained in the kernel of the cup product 
map.
\end{lem}
\begin{proof}
We prove this by direct computations.

\bigskip
\noindent
$[42]:$ \ Consider $v_{[2^2]}=T^H (a_1,a_2,a_1,a_2) \in [2^2]$ and 
$\frac{1}{2}V_2 (v_{[2^2]})=T^H (a_1,b_2,a_1,a_2) \in [2^2]$. 
Then 
\[[T^H (a_1,a_2,a_1,a_2),T^H (a_1,b_2,a_1,a_2)]=
2T(a_1,a_2,a_1,a_1,a_2,a_1),\]
where we put 

\vspace*{15pt}
\hspace*{30pt}
%WinTpicVersion3.08
\unitlength 0.1in
\begin{picture}( 38.9500,  5.8500)(  0.0500,-10.0500)
% STR 2 0 3 0
% 3 1800 900 1800 1000 5 0
% $a$
\put(18.0000,-10.0000){\makebox(0,0){$a$}}%
% STR 2 0 3 0
% 3 2300 405 2300 505 5 0
% $b$
\put(23.0000,-5.0500){\makebox(0,0){$b$}}%
% STR 2 0 3 0
% 3 2700 405 2700 505 5 0
% $c$
\put(27.0000,-5.0500){\makebox(0,0){$c$}}%
% STR 2 0 3 0
% 3 3100 405 3100 505 5 0
% $d$
\put(31.0000,-5.0500){\makebox(0,0){$d$}}%
% STR 2 0 3 0
% 3 950 700 950 800 5 0
% $T(a,b,c,d,e,f) \ :=$
\put(9.5000,-8.0000){\makebox(0,0){$T(a,b,c,d,e,f) \ :=$}}%
% STR 2 0 3 0
% 3 4600 700 4600 800 5 0
% $\in \ \DHgb^\Q (4).$
\put(46.0000,-8.0000){\makebox(0,0){$\in \ \DHgb^\Q (4).$}}%
% LINE 1 0 3 0
% 10 1900 1005 3900 1005 3500 1005 3500 605 
%3100 605 3100 1005 2700 1005 2700 605 2300 605 2300 1005
% 
\special{pn 13}%
\special{pa 1900 1006}%
\special{pa 3900 1006}%
\special{fp}%
\special{pa 3500 1006}%
\special{pa 3500 606}%
\special{fp}%
\special{pa 3100 606}%
\special{pa 3100 1006}%
\special{fp}%
\special{pa 2700 1006}%
\special{pa 2700 606}%
\special{fp}%
\special{pa 2300 606}%
\special{pa 2300 1006}%
\special{fp}%
% STR 2 0 3 0
% 3 3500 405 3500 505 5 0
% $e$
\put(35.0000,-5.0500){\makebox(0,0){$e$}}%
% STR 2 0 3 0
% 3 4000 900 4000 1000 5 0
% $f$
\put(40.0000,-10.0000){\makebox(0,0){$f$}}%
\end{picture}%

\vspace*{15pt}
\noindent
Then we have 
\[p_6^{(1,2)(3,4)(5)(6)} \circ \iota_6 
(2T(a_1,a_2,a_1,a_1,a_2,a_1))=-60 (a_1 \wedge a_2) \otimes 
(a_1 \wedge a_2) \otimes a_1 \otimes a_1.\] 
This is the highest weight vector of $[42]$, so that $[42]$ is contained 
in $\Im [\cdot,\cdot]$ for $g \ge 2$.

\bigskip
\noindent
$[31^3]:$ \ Consider $T^H (a_1,b_2,a_1,a_2), \ 
T^H (a_1,a_2,a_3,a_4) \in [2^2]$. Then 
\[[T^H (a_1,b_2,a_1,a_2), 
T^H (a_1,a_2,a_3,a_4)] =
T(a_1,a_2,a_1,a_1,a_4,a_3)\]
Applying $p_6^{(1,2,3,4)(5)(6)} \circ \iota_6$, we have 
\[-12 (a_1 \wedge a_2 \wedge a_3 \wedge a_4) \otimes a_1 \otimes a_1.\]
This is the highest weight vector of $[3 1^3]$, so that $[3 1^3]$ 
is contained in $\Im [\cdot,\cdot]$ for $g \ge 4$.

\bigskip
\noindent
$[2^3]:$ \ Consider $T^H (a_3,a_2,a_3,a_2) \in [2^2]$. Then 
\[[T^H (a_3,a_2,a_3,a_2),T^H (a_1,b_2,a_1,a_2)] =
2T(a_3,a_2,a_3,a_1,a_2,a_1)\]
Applying $p_6^{(1,2,3)(4,5,6)}\circ \iota_6$, we have 
\[-72 (a_1 \wedge a_2 \wedge a_3) \otimes 
(a_1 \wedge a_2 \wedge a_3).\]
This is the highest weight vector of $[2^3]$, so that $[2^3]$ is contained 
in $\Im [\cdot,\cdot]$ for $g \ge 3$.

\bigskip
\noindent
$[31]:$ \ Consider the vectors $T^H (a_1,a_2,b_3,a_2),\ 
T^H (a_1,b_2,a_1,a_3),\ T^H (a_3,a_2,a_3,a_2) \in [2^2]$ and 
$T^H (a_1,b_1,a_1,b_2) \in [2^2]+[1^2]$. Then we have
\begin{align*}
\xi_1 :&= [T^H (a_1,a_2,b_3,a_2),T^H (a_1,b_2,a_1,a_3)] \\
&=T(a_1,a_2,b_3,a_1,a_3,a_1)+T(b_3,a_2,a_1,a_1,a_3,a_1)
+T(a_1,a_2,a_2,a_1,b_2,a_1), \\
\xi_2 :&= \frac{1}{2}[T^H (a_1,a_2,a_1,a_2),T^H (a_1,b_1,a_1,b_2)] \\
&=T(a_2,a_1,a_2,a_1,b_2,a_1)+T(a_1,b_1,a_1,a_1,a_2,a_1)
\end{align*}
and
\[\begin{array}{|c||c|c|}
\hline
& \xi_1 & \xi_2 \\
\hline \hline
p_4^{(1,2)(3)(4)} \circ C_6^{(1,2)} \circ \iota_6 & 
0 & -12 (a_1 \wedge a_2) \otimes a_1 \otimes a_1 \\
\hline 
p_4^{(1,2)(3)(4)} \circ C_6^{(1,3)} \circ \iota_6 & 
4 (a_1 \wedge a_2) \otimes a_1 \otimes a_1 
& -4 (a_1 \wedge a_2) \otimes a_1 \otimes a_1 \\
\hline
\end{array}\]
Note that $\xi_1$ are defined for $g \ge 3$ and $\xi_2$ are for $g \ge 2$. 
From the table, we see that $\Im [\cdot,\cdot]$ contains $[31]$ 
for $g =2$, and $2[31]$ for $g \ge 3$.

\bigskip
\noindent
$[2]:$ \ Consider the vectors $\frac{1}{4}V_1V_2^2(v_{[2^2]})
=T^H (a_1,b_2,b_1,b_2) \in [2^2]$ and 
$T^H (a_1,b_1,a_1,b_2)$, $T^H (a_1,a_2,a_1,b_1) \in [2^2]+[1^2]$. 
Then we have
\begin{align*}
\xi_1 :&= \frac{1}{2}[T^H (a_1,a_2,a_1,a_2),T^H (a_1,b_2,b_1,b_2)] \\
&=T(a_1,b_2,b_2,a_2,a_2,a_1)+T(a_1,a_2,a_1,a_1,b_2,b_1)
+T(a_1,b_2,b_1,a_1,a_2,a_1), \\
\xi_2 :&= [T^H (a_1,b_1,a_1,b_2), T^H (a_1,a_2,a_1,b_1)] \\
&=T(b_2,a_1,b_1,a_1,a_2,a_1)+T(b_1,a_1,b_2,a_1,a_2,a_1)
+T(a_1,a_2,b_1,a_1,b_2,a_1) \\
&\quad +T(a_1,b_1,a_2,a_1,b_2,a_1)+T(b_1,a_1,a_1,a_1,b_1,a_1)
\end{align*}
and
\[\begin{array}{|c||c|c|}
\hline
& \xi_1 & \xi_2 \\
\hline \hline
C_4^{(1,2)} \circ C_6^{(1,2)} \circ \iota_6 & 
0 & 18 a_1 \otimes a_1 \\
\hline 
C_4^{(1,3)} \circ C_6^{(1,2)} \circ \iota_6 & 
6 a_1 \otimes a_1 
& -6 a_1 \otimes a_1 \\
\hline 
C_4^{(1,2)} \circ C_6^{(1,3)} \circ \iota_6 & 
-6 a_1 \otimes a_1 
& 6 a_1 \otimes a_1 \\
\hline
\end{array}\]
Note that $\xi_1,\xi_2$ are defined for $g \ge 2$. 
From the table, we see that $\Im [\cdot,\cdot]$ contains 
$2[2]$ for $g \ge 2$.

\bigskip
\noindent
$[21^2]:$ \ Consider $\Phi_2 (a_3 \wedge b_2)=
\sum_{i=1}^g T^H (a_3,b_2,a_i,b_i) \in [1^2]$.  
Then we have
\begin{align*}
&\frac{1}{2}\left[T^H (a_1,a_2,a_1,a_2),\ 
\sum_{i=1}^g T^H (a_3,b_2,a_i,b_i)\right] \\
&=\sum_{i=1}^gT(a_1,a_2,a_1,a_3,b_i,a_i)+T(a_3,b_2,a_2,a_1,a_2,a_1)
+T(a_2,a_1,a_2,a_1,b_2,a_3).
\end{align*}
Applying $p_4^{(1,2,3)(4)} \circ C_6^{(1,2)} \circ \iota_6$, we have 
\[6g (a_1 \wedge a_2 \wedge a_3) \otimes a_1.\]
This is the highest weight vector of $[21^2]$, so that $[21^2]$ is contained 
in $\Im [\cdot,\cdot]$ for $g \ge 3$.
\end{proof}

\subsection{An upper bound of the kernel}
\label{subsec:upperSg1}
We can detect summands in 
$\Im(\cup : \wedge^2 ([2^2]+[1^2]+[0]) \to H^2(\Kgb;\Q))$ by 
observing the dual 
\[\tau_{g,1}(2)_\ast : H_2 (\Kgb;\Q) 
\longrightarrow \wedge^2 ([2^2]+[1^2]+[0]) 
=H_2 (\DHgb (2);\Q).\]

For each pair $(\varphi, \psi)$ of elements of $\Kgb$ which are 
commutative, we have 
a homomorphism $f:\Z^2 \to \Kgb$ sending the standard generators 
of $\Z^2$ to $\varphi$ and $\psi$. Then we can construct an element of 
$H_2 (\Kgb;\Q)$ by considering the image of the generator of 
$1 \in \Q \cong H_2 (\Z^2;\Q)$ by $f_\ast$. 
Such a class is called an {\it abelian cycle}. The following is 
easily proved (see \cite[Section 2]{sa}). 
\begin{lem}\label{abelian}
Let $A$ be a finitely generated free abelian group and $f:\Z^2 \to A$ 
be a group homomorphism. Then the image of the generator 
$1 \in \Q \cong H_2 (\Z^2;\Q)$ by $f_\ast$ is 
\[f(e_1) \wedge f(e_2) \in \wedge^2 (A \otimes \Q) 
\cong H_2 (A;\Q),\]
where $e_1$, $e_2$ are the standard generators of $\Z^2$
\end{lem}
In our computation, we use this lemma for the image of 
each abelian cycle of $\Kgb$ by $\tau_{g,1}(2)_\ast$. 
\begin{lem}\label{lem:upperSg1}
For $g \ge 4$, $\Im (\tau_{g,1}(2)_\ast : H_2 (\Kgb;\Q) \to 
H_2 (\DHgb (2);\Q))$ contains 
\[[431]+[32^21]+2[321]+2[21^2]+[3^2]+[2^21^2]+2[2^2]+2[1^2].\]
Hence, these summands are not contained in the kernel of the cup product 
map.
\end{lem}
\begin{proof}
We first show that the following elements 
\begin{align*}
w_1&:=T^H(a_1,b_1,a_1,b_1) \wedge T^H(a_2,b_2,a_2,b_2)\\
w_2&:=T^H(a_1,b_1,a_1,b_1) \wedge T^H(a_3,b_3,a_3,b_3)\\
w_3&:=T^H(a_1,b_1,a_1,b_1) \wedge T^H(a_1,b_1,a_2,b_2)
\end{align*}
are contained in 
$\Im (\tau_{g,1}(2)_\ast : H_2 (\Kgb;\Q) \to 
H_2 (\DHgb (2);\Q))$. Indeed, 
the first two are easily obtained from Proposition \ref{image}. 
As for the last one, 
consider the abelian cycle corresponding to the pair of the 
elements $\psi_1$, $\psi_2 \in \Kgb$ in Section \ref{subsec:Johnson}. 
It is mapped by $\tau_{g,1}(2)_\ast$ to
\begin{align*}
\frac{1}{4}T^H&(a_1,b_1,a_1,b_1) \wedge 
\sum_{i,j=1}^2 T^H(a_i,b_i,a_j,b_j) \\
&=\frac{1}{2} T^H(a_1,b_1,a_1,b_1) \wedge T^H(a_1,b_1,a_2,b_2) 
+\frac{1}{4}T^H(a_1,b_1,a_1,b_1) \wedge T^H(a_2,b_2,a_2,b_2),
\end{align*}
and our claim follows from this. 

We decompose $w_1$, $w_2$, $w_3$ to each summands. 
We define a map
\[\iota :\wedge^2 \DHgb^\Q (2) \stackrel{\wedge^2 \iota_4}{\hookrightarrow} 
\wedge^2 (\Hq^{\otimes 4}) 
\hookrightarrow \Hq^{\otimes 8}\]
where the second inclusion sends $X \wedge Y$ to $X \otimes Y - Y \otimes X$ 
for $X,Y \in \Hq^{\otimes 4}$.

\bigskip
\noindent
$[431]:$ \ Consider $w_1=T^H(a_1,b_1,a_1,b_1) \wedge T^H(a_2,b_2,a_2,b_2)$. 
Then we have
\begin{align*}
w_1&\stackrel{X_{2,3}}{\longmapsto} 
-2 T^H(a_1,b_1,a_1,b_1) \wedge T^H(a_2,b_2,a_2,b_2) \\
&\stackrel{X_{1,2}^4}{\longmapsto} 
-48 T^H(a_1,b_2,a_1,b_2) \wedge T^H(a_1,b_3,a_1,b_2) \\
&\stackrel{U_2^3}{\longmapsto} 
-288 T^H(a_1,a_2,a_1,a_2) \wedge T^H(a_1,b_3,a_1,a_2) \\
&\stackrel{U_3}{\longmapsto} 
-288 T^H(a_1,a_2,a_1,a_2) \wedge T^H(a_1,a_3,a_1,a_2). 
\end{align*}
By $p_8^{(1,2,5)(4,3)(6,7)(8)} \circ \iota$, 
the last term is mapped to 
\[20736 (a_1 \wedge a_2 \wedge a_3) \otimes (a_1 \wedge a_2) 
\otimes (a_1 \wedge a_2) \otimes a_1,\]
which is the highest weight vector of $[431]$, 
so that $[431]$ is contained in $\Im \tau_{g,1}(2)_\ast$ for $g \ge 3$. 

\bigskip
\noindent
$[32^21]:$ \ 
We have
\begin{align*}
&p_8^{(1,2,5,6)(3,4,7)(8)} \circ \iota 
(U_4 X_{2,4} U_2 U_3^2 X_{1,2} X_{1,3}^2 w_1)\\
=\ & p_8^{(1,2,5,6)(3,4,7)(8)} \circ \iota 
(-8T^H(a_1,a_3,a_1,a_3) \wedge T^H(a_1,a_2,a_2,a_4)) \\
=\ & 576 (a_1 \wedge a_2 \wedge a_3 \wedge a_4) \otimes 
(a_1 \wedge a_2 \wedge a_3) \otimes a_1,
\end{align*}
which is the highest weight vector of $[32^21]$, 
so that $[32^21]$ is contained in $\Im \tau_{g,1}(2)_\ast$ for $g \ge 4$.

\bigskip
\noindent
$[321]:$ \ Consider $w_1=T^H(a_1,b_1,a_1,b_1) \wedge T^H(a_2,b_2,a_2,b_2)$ and 
$w_3=T^H(a_1,b_1,a_1,b_1) \wedge T^H(a_1,b_1,a_2,b_2)$. 
Then we have
\begin{align*}
w_1&\stackrel{U_2 U_3^2 X_{1,2} X_{1,3}^2}{\longmapsto} 
8 T^H(a_1,a_3,a_1,a_3) \wedge T^H(a_1,a_2,a_2,b_2) \\
&\hspace{18pt} \longmapsto \hspace{17pt} 
8 T^H(a_1,a_2,a_1,a_2) \wedge T^H(a_1,a_3,a_3,b_3) \ =: \xi_1, 
\end{align*}
where the second map is applying the element of $\Symp{Q}$ 
which acts on $\Hq$ by 
\[a_i \mapsto \left\{ 
\begin{array}{ll} a_3 & (i=2) \\ a_2 & (i=3) \\ a_i & (i \neq 2,3) 
\end{array}\right. , \qquad b_i \mapsto \left\{ 
\begin{array}{ll} b_3 & (i=2) \\ b_2 & (i=3) \\ b_i & (i \neq 2,3) 
\end{array}\right. , \]
and we denote it by $(a_2 \leftrightarrow a_3, b_2 \leftrightarrow b_3) 
\in \Symp{Q}$, for short. We also have
{\small 
\begin{align*}
w_3 &\ \stackrel{U_2 U_3^2 X_{1,2} X_{1,3}^2}{\longmapsto} 
\hspace{-10pt}
-4 T^H(a_1,a_3,a_1,a_3) \wedge T^H(a_1,a_2,a_2,b_2) 
+4 T^H(a_1,a_3,a_1,a_3) \wedge T^H(a_1,b_1,a_1,a_2) \\
&\hspace{45pt} -8 T^H(a_1,a_3,a_1,a_2) \wedge T^H(a_1,a_3,a_2,b_2) 
+8 T^H(a_1,a_3,a_1,b_1) \wedge T^H(a_1,a_3,a_1,a_2) \\
&\stackrel{(a_2 \leftrightarrow a_3, b_2 \leftrightarrow b_3)}{\longmapsto} 
\hspace{-10pt}
-4 T^H(a_1,a_2,a_1,a_2) \wedge T^H(a_1,a_3,a_3,b_3) 
+4 T^H(a_1,a_2,a_1,a_2) \wedge T^H(a_1,b_1,a_1,a_3) \\
&\hspace{45pt} -8 T^H(a_1,a_2,a_1,a_3) \wedge T^H(a_1,a_2,a_3,b_3) 
+8 T^H(a_1,a_2,a_1,b_1) \wedge T^H(a_1,a_2,a_1,a_3) \\
&\qquad =: \xi_2, 
\end{align*}}
and
{\small 
\[\begin{array}{|c||c|c|}
\hline
& \xi_1 & \xi_2 \\
\hline \hline
p_6^{(1,2,3)(4,5)(6)} \circ C_8^{(1,2)} \circ \iota & 
288 (a_1 \wedge a_2 \wedge a_3) \otimes (a_1 \wedge a_2) \otimes a_1 
& 240 (a_1 \wedge a_2 \wedge a_3) \otimes (a_1 \wedge a_2) \otimes a_1 \\
\hline 
p_6^{(1,2,4)(3,5)(6)} \circ C_8^{(1,7)} \circ \iota & 
0 & 120 (a_1 \wedge a_2 \wedge a_3) \otimes (a_1 \wedge a_2) \otimes a_1 \\
\hline
\end{array}\]
}
Note that $\xi_1, \xi_2$ are defined for $g \ge 3$. 
From the table, we see that $\Im \tau_{g,1}(2)_\ast$ 
contains $2[321]$ for $g \ge 3$.

\bigskip
\noindent
$[21^2]:$ \ Consider $w_1$ and $w_3$. Then we have
{\small 
\begin{align*}
w_1&\stackrel{U_2 U_3 X_{1,2} X_{1,3}}{\longmapsto} 
2 T^H(a_1,a_2,a_1,a_3) \wedge T^H(a_2,b_2,a_2,b_2) 
-4 T^H(a_1,b_1,a_1,a_3) \wedge T^H(a_1,a_2,a_2,b_2) \\
&\qquad =: \xi_1, \\
w_3 &\stackrel{U_2 U_3 X_{1,2} X_{1,3}}{\longmapsto} 
2 T^H(a_1,a_2,a_1,a_3) \wedge T^H(a_1,b_1,a_2,b_2) 
+2 T^H(a_1,b_1,a_1,a_3) \wedge T^H(a_1,a_2,a_2,b_2) \\
&\hspace{43pt} -2 T^H(a_1,b_1,a_1,a_3) \wedge T^H(a_1,b_1,a_1,a_2) 
+2 T^H(a_1,b_1,a_1,a_2) \wedge T^H(a_1,a_3,a_2,b_2) \\
&\hspace{43pt} -T^H(a_1,b_1,a_1,b_1) \wedge T^H(a_1,a_3,a_1,a_2) \\
&\qquad =: \xi_2, 
\end{align*}}
and
\[\begin{array}{|c||c|c|}
\hline
& \xi_1 & \xi_2 \\
\hline \hline
p_4^{(1,2,3)(4)} \circ C_6^{(5,6)} \circ C_8^{(1,2)} \circ \iota & 
-144 (a_1 \wedge a_2 \wedge a_3) \otimes a_1 
& -48 (a_1 \wedge a_2 \wedge a_3) \otimes a_1 \\
\hline 
p_4^{(1,2,3)(4)} \circ C_6^{(3,5)} \circ C_8^{(1,7)} \circ \iota & 
0 & 24 (a_1 \wedge a_2 \wedge a_3) \otimes a_1 \\
\hline
\end{array}\]
Note that $\xi_1, \xi_2$ are defined for $g \ge 3$. 
From the table, we see that $\Im \tau_{g,1}(2)_\ast$ 
contains $2[21^2]$ for $g \ge 3$.

\bigskip
\noindent
$[3^2]:$ \ We have 
{\small 
\[w_1 \stackrel{U_2^3 X_{1,2}^3}{\longmapsto} 
-72 T^H(a_1,a_2,a_1,b_1) \wedge T^H(a_1,a_2,a_1,a_2) 
+72 T^H(a_1,a_2,a_1,a_2) \wedge T^H(a_1,a_2,a_2,b_2) \]
}
Applying $p_6^{(1,2)(3,4)(5,6)} \circ C_8^{(1,2)} 
\circ \iota$, we obtain
\[-10368 (a_1 \wedge a_2) \otimes (a_1 \wedge a_2) 
\otimes (a_1 \wedge a_2).\]
Hence $\Im \tau_{g,1}(2)_\ast$ contains $[3^2]$ for $g \ge 2$.

\bigskip
\noindent
$[2^21^2]:$ \ We have 
\begin{align*}
w_1 &\stackrel{U_3 X_{1,3} U_4 X_{2,4} U_3 X_{1,3}}{\longmapsto} 
-4 T^H(a_1,a_3,a_1,a_3) \wedge T^H(a_2,b_2,a_2,a_4) \\
&\hspace{10pt}
\stackrel{(a_2 \leftrightarrow a_3,b_2 \leftrightarrow b_3)}
{\longmapsto} \hspace{10pt}
-4 T^H(a_1,a_2,a_1,a_2) \wedge T^H(a_3,b_3,a_3,a_4) \\
&\stackrel{p_6^{(1,2,3,4)(5,6)} \circ C_8^{(1,2)} 
\circ \iota}{\longmapsto} 
288 (a_1 \wedge a_2 \wedge a_3 \wedge a_4) \otimes (a_1 \wedge a_2). 
\end{align*}
Hence $\Im \tau_{g,1}(2)_\ast$ contains $[2^2 1^2]$ for $g \ge 4$.

\bigskip
\noindent
$[2^2]:$ \ We have
{\small 
\begin{align*}
w_2 &\stackrel{U_2^2 X_{1,2}^2}{\longmapsto} 
4 T^H(a_1,a_2,a_1,a_2) \wedge T^H(a_3,b_3,a_3,b_3) =: \xi_1, \\
w_3 &\stackrel{U_2^2 X_{1,2}^2}{\longmapsto} 
4 T^H(a_1,a_2,a_1,a_2) \wedge T^H(a_1,b_1,a_2,b_2) 
+8 T^H(a_1,b_1,a_1,a_2) \wedge T^H(a_1,a_2,a_2,b_2) \\
&\qquad -4 T^H(a_1,b_1,a_1,b_1) \wedge T^H(a_1,a_2,a_1,a_2) \\
&\quad =: \xi_2, 
\end{align*}
}
and
\[\begin{array}{|c||c|c|}
\hline
& \xi_1 & \xi_2 \\
\hline \hline
p_4^{(1,2)(3,4)} \circ C_6^{(1,2)} \circ C_8^{(1,2)} \circ \iota & 
-576 (a_1 \wedge a_2) \otimes (a_1 \wedge a_2)
& -960 (a_1 \wedge a_2) \otimes (a_1 \wedge a_2) \\
\hline 
p_4^{(1,2)(3,4)} \circ C_6^{(1,2)} \circ C_8^{(1,6)} \circ \iota & 
0 & -240 (a_1 \wedge a_2) \otimes (a_1 \wedge a_2) \\
\hline
\end{array}\]
Note that $\xi_1$ are defined for $g \ge 3$, and $\xi_2$ for $g \ge 2$. 
From the table, we see that $\Im \tau_{g,1}(2)_\ast$ contains 
$[2^2]$ for $g = 2$, and $2[2^2]$ for $g \ge 3$.

\bigskip
\noindent
$[1^2]:$ \ We have
{\small 
\begin{align*}
w_1 &\stackrel{U_2 X_{1,2}}{\longmapsto} 
-2 T^H(a_1,b_1,a_1,a_2) \wedge T^H(a_2,b_2,a_2,b_2) 
+2 T^H(a_1,b_1,a_1,b_1) \wedge T^H(a_1,a_2,a_2,b_2) \\
&\quad =: \xi_1, \\
w_3 &\stackrel{U_2 X_{1,2}}{\longmapsto} 
-2 T^H(a_1,b_1,a_1,a_2) \wedge T^H(a_1,b_1,a_2,b_2) 
-T^H(a_1,b_1,a_1,b_1) \wedge T^H(a_1,a_2,a_2,b_2) \\
&\hspace{28pt} + T^H(a_1,b_1,a_1,b_1) \wedge T^H(a_1,b_1,a_1,a_2) \\
&\quad =: \xi_2, 
\end{align*}
}
and
\[\begin{array}{|c||c|c|}
\hline
& \xi_1 & \xi_2 \\
\hline \hline
p_2^{(1,2)} \circ C_4^{(1,2)} \circ C_6^{(1,2)} 
\circ C_8^{(1,2)} \circ \iota & 
288 a_1 \wedge a_2 & 96 a_1 \wedge a_2 \\
\hline 
p_2^{(1,2)} \circ C_4^{(1,2)} \circ C_6^{(1,2)} 
\circ C_8^{(1,5)} \circ \iota & 
0 & -48 a_1 \wedge a_2 \\
\hline
\end{array}\]
Note that $\xi_1, \xi_2$ are defined for $g \ge 3$. 
From the table, we see that $\Im \tau_{g,1}(2)_\ast$ contains 
$2[1^2]$ for $g \ge 2$.
\end{proof}
From the above arguments, Theorem \ref{mainSg1} follows.

\subsection{Statement for $\Kgp$}\label{subsec:compSgp}
As seen in Section \ref{subsec:Johnson-and-rep}, 
we have $\Im \tau_{g,\ast}^\Q (2) = \DHgp^\Q (2) \cong [2^2]+[1^2]$, 
so that an injection $[2^2]+[1^2] \hookrightarrow H^1 (\Kgp;\Q)$ 
is obtained. We now consider the cup product map
\[\cup : \wedge^2 ([2^2]+[1^2]) \longrightarrow H^2 (\Kgp;\Q).\]
Note that 
\begin{align*}
\wedge^2 ([2^2]+[1^2]) 
\cong& \wedge^2 [2^2] + [2^2] \otimes [1^2] + \wedge^2 [1^2] \\
\cong&\ ([431]+[42]+[3 2^2 1]+[321]+[31^3]+
[31]+[2^3]+[21^2]+[2]) \\
&+ ([321]+[31]+[21^2]+[3^2]+[2^2 1^2]+[2^2]+[1^2]) + ([21^2]+[2]).
\end{align*}
\begin{thm}\label{mainSgp}
For $g \ge 4$, the kernel of the cup product map 
$\cup : \wedge^2 ([2^2]+[1^2]) \to H^2(\Kgp;\Q)$ is 
\[[42]+[31^3]+2[31]+[2^3]+[21^2]+2[2], \]
which is, as an $\Symp{Q}$-vector space, 
isomorphic to $\Im \tau_{g,1}^\Q (4) \cong \Im\tau_{g,\ast}^\Q (4)$.
\end{thm}
\begin{proof}
Since $\Ker(\wedge^2 \DHgb^\Q (2) \to \wedge^2 \DHgp^\Q (2)) = 
1[2^2]+1[1^2]$, only we have to do is 
to observe the summands $[2^2]$ and $[1^2]$. 
In the last subsection, we have proved that $\Im \tau_{g,1}(2)_\ast$ 
contains $2[2^2]+2[1^2]$. Therefore 
$\Im \tau_{g,\ast}(2)_\ast$ certainly contains $1[2^2]+1[1^2]$, 
and this completes the proof. 
\end{proof}

\vspace{10pt}

\section{Computation for the case of a closed surface}\label{sec:compSg}
\subsection{Statement for $\Kg$}\label{subsec:statementSg}

As seen in Section \ref{subsec:Johnson-and-rep}, 
we have $\Im \tau_g^\Q (2) = \DHg^\Q (2) \cong [2^2]$, 
so that an injection $[2^2] \hookrightarrow H^1 (\Kg;\Q)$ 
is obtained. We now consider the cup product map
\[\cup : \wedge^2 [2^2] \longrightarrow H^2 (\Kg;\Q).\]
Note that 
$\wedge^2 [2^2] \cong 
[431]+[42]+[3 2^2 1]+[321]+[31^3]+[31]+[2^3]+[21^2]+[2]$.
\begin{thm}\label{mainSg}
For $g \ge 4$, the kernel of the cup product map 
$\cup : \wedge^2 [2^2] \to H^2(\Kg;\Q)$ is 
\[[42]+[31^3]+[31]+[2^3]+[2], \]
which is, as an $\Symp{Q}$-vector space, 
isomorphic to $\Im \tau_g^\Q (4)$.
\end{thm}
This theorem follows from the arguments in the next two subsections. 

\subsection{A lower bound of the kernel}
\label{subsec:lowerSg}
By an argument similar to that in Section \ref{subsec:lowerSg1}, 
we can find some summands in 
$\wedge^2 [2^2]$ which 
vanish in $H^2 (\Kg;\Q)$ by showing the following. 
\begin{lem}\label{lem:lowerSg}
For $g \ge 4$, the image of the Lie bracket 
$[\cdot,\cdot] : \wedge^2 \Im \tau_g^\Q (2) \to 
\Im \tau_g^\Q (4)$ contains 
\[[42]+[31^3]+[31]+[2^3]+[2].\]
\end{lem}
\begin{proof}
In the previous section, we have proved that
\[\Im([\cdot,\cdot] :\wedge^2 \Im \tau_{g,1}^\Q (2) 
\to \Im \tau_{g,1}^\Q (4))=
[42]+[31^3]+2[31]+[2^3]+[21^2]+2[2] \subset \DHgb^\Q (4).\]
On the other hand, $\Ker(\DHgb^\Q (4) \to \DHg^\Q (4))=[31]+[21^2]+2[2]$. 
Hence at least $[42]+[31^3]+[31]+[2^3]$ is contained in 
$\Im[\cdot,\cdot] \subset \Im \tau_{g}^\Q (4)$. 
We now show that $[2]$ is also contained in $\Im[\cdot,\cdot]$. 

The multiplicity of $[2]$ in $\DHgb^\Q (4)$ is 3. We now prepare the 
following three vectors in $\DHgb^\Q (4)$:
\begin{align*}
&\chi_1:=\frac{1}{2}[T^H (a_1,a_2,a_1,a_2),T^H (a_1,b_2,b_1,b_2)],\\
&\chi_2:=\frac{1}{g-1}\Phi_4 
\left(\sum_{i=1}^g \ [[a_i,a_1],[b_i,a_1]] \right),\\
&\chi_3:=\frac{1}{g-1}[\theta_1,\theta_2]
\end{align*}
where

\vspace*{15pt}
\hspace*{-80pt}
%WinTpicVersion3.08
\unitlength 0.1in
\begin{picture}( 57.0500,  9.0000)( -7.0500,-13.2000)
% STR 2 0 3 0
% 3 1050 800 1050 900 5 0
% $\theta_1 := \displaystyle\sum_{i=1}^g$
\put(10.5000,-9.0000){\makebox(0,0){$\theta_1 := 
\displaystyle\sum_{i=1}^g$}}%
% LINE 1 0 3 0
% 2 2000 805 2000 605
% 
\special{pn 13}%
\special{pa 2000 806}%
\special{pa 2000 606}%
\special{fp}%
% LINE 1 0 3 0
% 8 2000 805 1500 1305 1900 1305 1700 1105 2100 1305 
%2300 1105 2500 1305 2000 805
% 
\special{pn 13}%
\special{pa 2000 806}%
\special{pa 1500 1306}%
\special{fp}%
\special{pa 1900 1306}%
\special{pa 1700 1106}%
\special{fp}%
\special{pa 2100 1306}%
\special{pa 2300 1106}%
\special{fp}%
\special{pa 2500 1306}%
\special{pa 2000 806}%
\special{fp}%
% STR 2 0 3 0
% 3 2000 405 2000 505 5 0
% $a_1$
\put(20.0000,-5.0500){\makebox(0,0){$a_1$}}%
% STR 2 0 3 0
% 3 1500 1305 1500 1405 5 0
% $a_i$
\put(15.0000,-14.0500){\makebox(0,0){$a_i$}}%
% STR 2 0 3 0
% 3 1900 1305 1900 1405 5 0
% $a_1$
\put(19.0000,-14.0500){\makebox(0,0){$a_1$}}%
% STR 2 0 3 0
% 3 2100 1305 2100 1405 5 0
% $b_i$
\put(21.0000,-14.0500){\makebox(0,0){$b_i$}}%
% STR 2 0 3 0
% 3 2500 1305 2500 1405 5 0
% $a_1$
\put(25.0000,-14.0500){\makebox(0,0){$a_1$}}%
% LINE 1 0 3 0
% 6 4700 900 4700 600 4700 900 4400 1200 5000 1200 4700 900
% 
\special{pn 13}%
\special{pa 4700 900}%
\special{pa 4700 600}%
\special{fp}%
\special{pa 4700 900}%
\special{pa 4400 1200}%
\special{fp}%
\special{pa 5000 1200}%
\special{pa 4700 900}%
\special{fp}%
% STR 2 0 3 0
% 3 5600 800 5600 900 5 0
% $\in \DHgb^\Q (1)$.
\put(56.0000,-9.0000){\makebox(0,0){$\in \DHgb^\Q (1)$.}}%
% STR 2 0 3 0
% 3 3900 800 3900 900 5 0
% $\theta_2 := \displaystyle\sum_{j=1}^g$
\put(39.0000,-9.0000){\makebox(0,0){$\theta_2 := 
\displaystyle\sum_{j=1}^g$}}%
% STR 2 0 3 0
% 3 5100 1210 5100 1310 5 0
% $b_j$
\put(51.0000,-13.1000){\makebox(0,0){$b_j$}}%
% STR 2 0 3 0
% 3 4300 1210 4300 1310 5 0
% $a_j$
\put(43.0000,-13.1000){\makebox(0,0){$a_j$}}%
% STR 2 0 3 0
% 3 4700 410 4700 510 5 0
% $b_1$
\put(47.0000,-5.1000){\makebox(0,0){$b_1$}}%
% STR 2 0 3 0
% 3 3000 800 3000 900 5 0
% $\in \DHgb^\Q (3)$,
\put(30.0000,-9.0000){\makebox(0,0){$\in \DHgb^\Q (3)$,}}%
\end{picture}%

\vspace{15pt}
\noindent
Note that $\chi_2, \chi_3 \in 
\Ker(\DHgb^\Q (4) \to \DHg^\Q (4))$. 
Indeed the commutative diagram 
\[\begin{CD}
@. \Lgb^\Q (4) @>\Phi_4>> \DHgb^\Q (4) @. @. \\
@. @VVV @VVV @. @. \\
0 @>>> \Lg^\Q (4) @>\Psi_4>> \DHgp^\Q (4) @>>> \DHg^\Q (4) @>>> 0,
\end{CD}\]
whose bottom row is exact, shows $\chi_2 \in 
\Ker(\DHgb^\Q (4) \to \DHg^\Q (4))$. Also, 
$\theta_2 \in [1] = \Ker(\DHgb^\Q (1) \to \DHg^\Q (1)) 
\subset \DHgb^\Q (1)$ and the fact 
that the map $\DHgb^\Q \to \DHg^\Q$ preserves brackets show $\chi_3 \in 
\Ker(\DHgb^\Q (4)$ $\to \DHg^\Q (4))$. 
By direct calculations, we have 
\[\begin{array}{|c||c|c|c|}
\hline
& \chi_1 & \chi_2 & \chi_3\\
\hline \hline
C_4^{(1,2)} \circ C_6^{(1,2)} \circ \iota_6 & 
0 & 2a_1 \otimes a_1 & (-8g-2) a_1 \otimes a_1\\
\hline 
C_4^{(1,3)} \circ C_6^{(1,2)} \circ \iota_6 & 
6 a_1 \otimes a_1 
& (4g-2)a_1 \otimes a_1 & (12g-2) a_1 \otimes a_1\\
\hline 
C_4^{(1,2)} \circ C_6^{(1,3)} \circ \iota_6 & 
-6 a_1 \otimes a_1 & -2 a_1 \otimes a_1
& -10 a_1 \otimes a_1 \\
\hline
\end{array}\]
and from this we observe that $\chi_1,\chi_2,\chi_3$ generate 
$3[2]$ in $\DHgb^\Q (4)$. 

Combining the above, we see that $\chi_1$ survives 
in $\Im \tau_{g}^\Q (4) \subset \DHg^\Q (4)$. 
\end{proof}

\subsection{An upper bound of the kernel}
\label{subsec:upperSg}
As in Section \ref{subsec:upperSg1}, 
we can find summands in 
$\wedge^2 [2^2]$ which 
survive in $H^2 (\Kg;\Q)$ by showing the following. 
\begin{lem}\label{lem:upperSg}
For $g \ge 4$, $\Im (\tau_g (2)_\ast : H_2 (\Kg;\Q) \to 
H_2 (\DHg (2);\Q))$ contains 
\[[431]+[32^21]+[321]+[21^2].\]
\end{lem}
\begin{proof}
In the previous section, we have proved that 
\[\Im(\tau_{g,1}(2)_\ast : H_2 (\Kgb;\Q) \to H_2 (\DHgb (2);\Q))
\supset [431]+[32^21]+2[321]+2[21^2].\]
On the other hand, we have
\begin{align*}
&\Ker(H_2 (\DHgb (2);\Q) \to H_2 (\DHg (2);\Q)) \\
=\ &\Ker(\wedge^2 \DHgb^\Q (2) \to \wedge^2 \DHg^\Q (2)) \\
=\ &[2^2] \otimes [1^2] +\wedge^2 [1^2]+ [2^2] \otimes [0] 
+[1^2] \otimes [0]\\
=\ &([321]+[31]+[21^2]+[3^2]+[2^2 1^2]+[2^2]+[1^2]) + 
([21^2]+[2]) + [2^2] + [1^2].
\end{align*}
Hence at least $[431]+[32^21]$ is contained in 
$\Im \tau_g (2)_\ast$. $[321]$ is 
also contained in $\Im\tau_g (2)_\ast$, since 
the multiplicity of $[321]$ in $\wedge^2 \DHgb^\Q (2)$ is 2 in which 
that in $\Ker(H_2 (\DHgb (2);\Q) \to H_2 (\DHg (2);\Q))$ is 1. 
We now show that $[21^2]$ is certainly contained in 
$\Im\tau_g (2)_\ast$. 

By Proposition \ref{image}, the projection 
$\wedge^2 \DHgb^\Q (2) \to \wedge^2 \DHg^\Q (2)$ maps 
$T^H (a_1,b_1,a_1,b_1) \wedge T^H (a_2,b_2,$ $a_2,b_2) \in 
\Im \tau_{g,1}(2)_\ast$ to 
\begin{align*}
&\left\{ T^H (a_1,b_1,a_1,b_1)
-\frac{3}{g+1}\Phi_2 (a_1 \wedge b_1)
+\frac{3}{(2g+1)(g+1)}\Phi_2 (\omega_0) \right\} \\
&\wedge \left\{ T^H (a_2,b_2,a_2,b_2)
-\frac{3}{g+1}\Phi_2 (a_2 \wedge b_2)
+\frac{3}{(2g+1)(g+1)}\Phi_2 (\omega_0) \right\} \\
\\
=\ & \ T^H (a_1,b_1,a_1,b_1) \wedge T^H (a_2,b_2,a_2,b_2) \\
&-\frac{3}{g+1}\left\{(T^H (a_1,b_1,a_1,b_1) \wedge \Phi_2 (a_2 \wedge b_2)
+\Phi_2 (a_1 \wedge b_1) \wedge T^H (a_2,b_2,a_2,b_2) \right\}\\
&+\frac{9}{(g+1)^2} \Phi_2 (a_1 \wedge b_1) \wedge \Phi_2 (a_2 \wedge b_2)\\
&+(\mbox{the other summands}),
\end{align*}
where the terms after the third are contained in 
$[2^2] \otimes [0]+[1^2] \otimes [0]$, and 
do not contribute to $[21^2]$. Then we have 
\begin{align*}
T^H (a_1,b_1,&a_1,b_1) \wedge T^H (a_2,b_2,a_2,b_2) \\
\stackrel{U_2 U_3 X_{1,2} X_{1,3}}{\longmapsto} \qquad 
&2 T^H(a_1,a_2,a_1,a_3) \wedge T^H(a_2,b_2,a_2,b_2) \\
&\hspace{-12pt} -4 T^H(a_1,b_1,a_1,a_3) \wedge T^H(a_1,a_2,a_2,b_2) \\
\stackrel{p_4^{(1,2,3)(4)} \circ C_6^{(3,5)} \circ C_8^{(1,7)} 
\circ \iota}{\longmapsto}
&\ 0,\\
\\
-\frac{3}{g+1}\{T^H &(a_1,b_1,a_1,b_1) \wedge \Phi_2 (a_2 \wedge b_2)
+\Phi_2 (a_1 \wedge b_1) \wedge T^H (a_2,b_2,a_2,b_2) \}\\
= \quad \qquad-&\frac{3}{g+1}\sum_{i=1}^g
\left\{\begin{array}{l}
T^H (a_1,b_1,a_1,b_1) \wedge T^H (a_2,b_2,a_i,b_i) \\
+T^H (a_1,b_1,a_i,b_i)\wedge T^H (a_2,b_2,a_2,b_2) 
\end{array}\right\}\\
\stackrel{U_2 U_3 X_{1,2} X_{1,3}}{\longmapsto} \quad
-&\frac{6}{g+1}\sum_{i=1}^g
\left\{\begin{array}{l}
T^H (a_1,a_2,a_1,a_3) \wedge T^H (a_2,b_2,a_i,b_i) \\
-T^H (a_1,b_1,a_1,a_3)\wedge T^H (a_1,a_2,a_i,b_i) \\
-T^H (a_1,a_3,a_i,b_i)\wedge T^H (a_1,a_2,a_2,b_2)
\end{array}\right\}\\
\stackrel{p_4^{(1,2,3)(4)} \circ C_6^{(3,5)} \circ C_8^{(1,7)} 
\circ \iota}{\longmapsto}
& -\frac{144}{g+1}(a_1 \wedge a_2 \wedge a_3) \otimes a_1
\intertext{and} 
\displaystyle\frac{9}{(g+1)^2} 
\Phi_2 (& a_1 \wedge b_1) \wedge \Phi_2 (a_2 \wedge b_2)
\stackrel{X_{1,2} X_{1,3}}{\longmapsto} 0.
\end{align*}
Therefore $[21^2]$ 
is contained in $\Im\tau_g (2)_\ast$ for $g \ge 3$. 
\end{proof}

\noindent
From the above arguments, Theorem \ref{mainSg} follows. 

\begin{remark}
Brendle-Farb \cite{bf} studied $H^2 (\Kgb;\Z)$ by using 
the Birman-Craggs-Johnson 
homomorphism and its integral lift 
by Morita \cite{mo9}, and showed that the rank of $H^2 (\Kgb;\Z)$ is 
at least $16g^4+O(g^3)$. From our computation, 
we can give a sharper estimate. Indeed the summand 
$[32^21] \subset \wedge^2 [2^2]$ survives in $H^2 (\Kgb;\Q)$, 
$H^2 (\Kgp;\Q)$ and $H^2 (\Kg;\Q)$. 
The dimension of this summand is 
\[\frac{1}{36}(g-3)(g-2)(g-1)(g+2)(2g-1)(2g+1)^2(2g+3).\]
However it seems that we should not compare these results by 
simply seeing their orders with respect to $g$, 
since the classes detected in \cite{bf} come 
from a different context from ours. 
\end{remark}

\vspace{10pt}

\section{The cases of $g=2,3$}\label{sec:g23}
We state the corresponding theorems for the cases of $g=2,3$. Note that 
we have already proved them in the argument of the previous sections. 
\subsection{The case of $g=2$}
\begin{lem}
For $g =2$, the irreducible decompositions 
of $\wedge^2 [2^2]$, $[2^2] \otimes [1^2]$, $\wedge^2 [1^2]$, 
$[2^2] \otimes [0]$ and $[1^2] \otimes [0]$ are 
given by the following table. 

\[\begin{array}{c}
\begin{array}{|c||c|c|c|c|c|}
\hline
 & \wedge^2 [2^2] &  [2^2] \otimes [1^2]  
 & \wedge^2 [1^2] &  [2^2] \otimes [0]  
 &  [1^2] \otimes [0]  \\
\hline
 [42]  &  1  & & & & \\
\hline
 [3 1]  & &  1  & & & \\
\hline
 [2]  &  1  & &  1  & & \\
\hline
 [3^2]  & &  1  & & & \\
\hline
 [2^2]  & & &  &  1  & \\
\hline
 [1^2]  & &  1  & & &  1  \\
\hline
\end{array}\vspace{2mm}\\
\mbox{\footnotesize {\rm Table 2}. Irreducible decompositions 
of $Sp(4,\Q)$-modules} 
\end{array}\]
\end{lem}

\smallskip
\begin{thm}
For $g =2$, we have 
\begin{align*}
&\Ker (\cup:\wedge^2 ([2^2]+[1^2]+[0]) \to H^2 (\mathcal{K}_{2,1};\Q)) 
=[42]+[31]+2[2], \\
&\Ker (\cup:\wedge^2 ([2^2]+[1^2]) \to H^2 (\mathcal{K}_{2,\ast};\Q)) 
=[42]+[31]+2[2], \\
&\Ker (\cup:\wedge^2 [2^2] \to H^2 (\mathcal{K}_2;\Q)) =[42]+[2].
\end{align*}
\end{thm}
\noindent
Since $\mathcal{K}_2$ is a free group \cite{mess}, 
the third equality is trivial. 

\subsection{The case of $g=3$}
\begin{lem}
For $g = 3$, the irreducible decompositions 
of $\wedge^2 [2^2]$, $[2^2] \otimes [1^2]$, $\wedge^2 [1^2]$, 
$[2^2] \otimes [0]$ and $[1^2] \otimes [0]$ are 
given by the following table. 

\[\begin{array}{c}
\begin{array}{|c||c|c|c|c|c|}
\hline
 & \wedge^2 [2^2] &  [2^2] \otimes [1^2]  
 & \wedge^2 [1^2] &  [2^2] \otimes [0]  
 &  [1^2] \otimes [0]  \\
\hline
 [431]  &  1  & & & & \\
\hline
 [42]  &  1  & & & & \\
\hline
 [3 2 1]  &  1  &  1  & & & \\
\hline
 [3 1]  &  1  &  1  & & & \\
\hline
 [2^3]  &  1  & & & & \\
\hline
 [2 1^2]  &  1  &  1  &  1  & & \\
\hline
 [2]  &  1  & &  1  & & \\
\hline
 [3^2]  & &  1  & & & \\
\hline
 [2^2]  & &  1  &  &  1  & \\
\hline
 [1^2]  & &  1  & & &  1  \\
\hline
\end{array}\vspace{2mm}\\
\mbox{\footnotesize {\rm Table 3}. Irreducible decompositions 
of $Sp(6,\Q)$-modules} \end{array}\]
\end{lem}

\smallskip
\begin{thm}
For $g =3$, we have 
\begin{align*}
&\Ker (\cup:\wedge^2 ([2^2]+[1^2]+[0]) \to H^2 (\mathcal{K}_{3,1};\Q))
=[42]+2[31]+[2^3]+[21^2]+2[2], \\
&\Ker (\cup:\wedge^2 ([2^2]+[1^2]) \to H^2 (\mathcal{K}_{3,\ast};\Q))
=[42]+2[31]+[2^3]+[21^2]+2[2], \\
&\Ker (\cup:\wedge^2 [2^2] \to H^2 (\mathcal{K}_3;\Q))
=[42]+[31]+[2^3]+[2].
\end{align*}
\end{thm}

\vspace{20pt}

\section{Appendix}\label{ch:appendix}
The following is a simple MATHEMATICA program which helps us to 
calculate the inclusion 
$\iota_k : \Lgb^\Q (k) \hookrightarrow \Hq^{\otimes k}$, 
the projections $p_k^\sigma$ and the contractions $C_k^{(i,j)}$. 

\medskip

{\scriptsize
\begin{verbatim}
tens[]:=1

tens[x___,c_ y_,z___]:=
  c*tens[x, y,z]/;NumberQ[c]

tens[x___,c_,z___]:=
  c*tens[x,z]/;NumberQ[c]

tens[x___,c_ +y_,z___]:=
  tens[x, c,z]+tens[x,y,z]

SetAttributes[tens,Flat]

brac[x_,y_]:=tens[x,y]-tens[y,x]

tree[a_,b_,c_,d_,e_,f_]:=
  Expand[tens[a,brac[brac[brac[brac[f,e],d],c],b]]-
      tens[b,brac[brac[brac[brac[f,e],d],c],a]]-
      tens[c,brac[brac[brac[f,e],d],brac[b,a]]]-
      tens[d,brac[brac[brac[b,a],c],brac[f,e]]]+
      tens[e,brac[brac[brac[brac[b,a],c],d],f]]-
      tens[f,brac[brac[brac[brac[b,a],c],d],e]]]

Ht[a_,b_,c_,d_]:=
  Expand[tens[a,brac[b,brac[c,d]]]+tens[b,brac[brac[c,d],a]]+
      tens[c,brac[d,brac[a,b]]]+tens[d,brac[brac[a,b],c]]]

genusNo[x_]:=
  If[StringTake[ToString[x],{1}]=="a",
    ToExpression[
      StringDrop[ToString[x],1]],-ToExpression[StringDrop[ToString[x],1]]]

contract[a_,b_]:=
  If[genusNo[a]+genusNo[b]\[Equal]0,
    If[genusNo[a]>0,1,-1],0]

basisRecover[x_]:=
  If[x>0,ToExpression[ToString[SequenceForm["a",x]]],
    ToExpression[ToString[SequenceForm["b",-x]]]]

Wedge[]:=1

Wedge[x___,c_ y_,z___]:=
  c*Wedge[x, y,z]/;NumberQ[c]

Wedge[x___,c_,z___]:=
  c*Wedge[x,z]/;NumberQ[c]

Wedge[x___,c_ +y_,z___]:=
  Wedge[x, c,z]+Wedge[x,y,z]

Wedge[x___,c_,y___,c_,z___]:=0

Wedge[x___,c_,y___,d_,z___]:=
  -Wedge[x,d,y,c,z]/;compare[genusNo[c],genusNo[d]]\[Equal]1

SetAttributes[Wedge,Flat]

SetAttributes[genusNo,Listable]

compare[x_,y_]:=
  Which[x\[Equal]1/3,0,
    y\[Equal]1/3,0,
    Abs[x]>Abs[y],1,
    Abs[x]<Abs[y],0,
    x<y,1]
\end{verbatim}
}

\medskip

In this program, \verb/tens/ means $\otimes$, 
\verb/Wedge/ means $\wedge$, the function \verb/tree[a_,b_,c_,d_,/
\verb/e_,f_]/ 
makes $T(a,b,c,d,e,f)$ in 
$\Hq^{\otimes 6}$, \verb/Ht[a_,b_,c_,d_]/ makes $T^H(a,b,c,d)$ 
in $\Hq^{\otimes 4}$. The function \verb/brac[x_,y_]/ is 
the bracket operation. We now give a sample computation using 
this program below. We calculate 
\[p_4^{(1,2)(3,4)} \circ C_6^{(1,2)}\circ C_8^{(1,2)} \circ \iota 
(4T^H(a_1,a_2,a_1,a_2) \wedge T^H(a_3,b_3,a_3,b_3))\]
of $[2^2]$ in Section \ref{subsec:upperSg1}. 
We first input the above program. 

\medskip

{\scriptsize
\begin{verbatim}
In[21]:=
xsi1=4brac[Ht[a1,a2,a1,a2],Ht[a3,b3,a3,b3]];

In[22]:=
ClearAttributes[tens,Flat];
xsi1/.tens[a_,b_,c_,d_,e_,f_,g_,h_]:>
    tens[contract[a,b],contract[c,d],Wedge[e,f],Wedge[g,h]]

SetAttributes[tens,Flat]

Out[23]=
-576 tens[a1\[Wedge]a2,a1\[Wedge]a2]
\end{verbatim}
}
\noindent
The output is $-576 (a_1 \wedge a_2) \otimes (a_1 \wedge a_2)$.

\vspace{10pt}

\section{Acknowledgement}\label{ch:acknowledge}
The author would like to express his gratitude 
to Professor Shigeyuki Morita 
for his encouragement and helpful suggestions as well as 
allowing the author to access his result of 
computations \cite{mop}. 
The author also would like to thank the referee for 
useful comments and suggestions. In particular, one of them 
corrects an insufficient point in the original argument given 
in Section \ref{subsec:lowerSg1}. (See Remark \ref{referee}.)

This research was supported by 
JSPS Research Fellowships for Young Scientists.

\end{document}